\documentclass{siamart250211}
\usepackage[utf8]{inputenc} 
\usepackage{newunicodechar} 
\usepackage{amssymb}
\usepackage{hyperref}
\usepackage{graphicx}
\usepackage{stmaryrd} 
\usepackage{bbold} 
\usepackage{xcolor}
\usepackage{bm} 

\usepackage{pgfplots} 
\usepackage{tikz}
\usetikzlibrary{3d,perspective}

\usepackage{easyReview} 

\newcommand{\cO}{\mathcal{O}}

\newcommand{\complex}{\mathbb{C}}

\newcommand{\cF}{\mathcal{F}}
\newcommand{\irange}[1]{\llbracket{#1}\rrbracket}
\newcommand{\hprod}{{.*}}

\newcommand{\vect}{\text{vec}\,}
\newcommand{\reshape}{\text{reshape}\,}
\newcommand{\diag}{\text{diag}\,}

\DeclareMathSymbol{\dvchar}{\mathord}{operators}{"3A}
\newcommand{\dv}{\mskip 1mu \dvchar \mskip 1mu}

\renewcommand{\mod}[2]{#1\, \text{\scriptsize{\%}}\, #2}

\newsiamremark{remark}{Remark}
\crefname{remark}{remark}{remarks}
\Crefname{remark}{Remark}{Remarks}

\title{Fast circular convolution with periodic kernels\thanks{Submitted to the editors
DATE.}}
\author{Fernando Guevara Vasquez\footnotemark[3]\ \thanks{Corresponding
Author. Department of Mathematics,
University of Utah, Salt Lake City, UT
(\email{fguevara@math.utah.edu})} \and Fernando González del
Cueto\thanks{Lumos Imaging Inc., 134 H St, Salt Lake City, UT}}

\headers{Fast convolution with periodic kernels}{F. González del Cueto and F. Guevara Vasquez}

\newlength{\imagesize} 
\begin{document}
\maketitle
\begin{abstract}
    We study the periodic circular convolution or PCC, where the
    convolution kernel depends periodically on the spatial location.  We
    introduce two fast algorithms to compute the PCC that rely on
    the Discrete Fourier Transform. One of the algorithms reveals that
    the matrix representation of the PCC is unitarily equivalent to a
    block diagonal matrix, which facilitates the analysis
    of PCCs and operations such as their product and inversion. We
    present two sample applications. One to deconvolve the action of an
    optical system whose point spread function depends periodically in
    space. The second is to compute dispersion relations for a
    periodic chain of masses.
\end{abstract}



\section{Introduction}
\label{sec:intro}
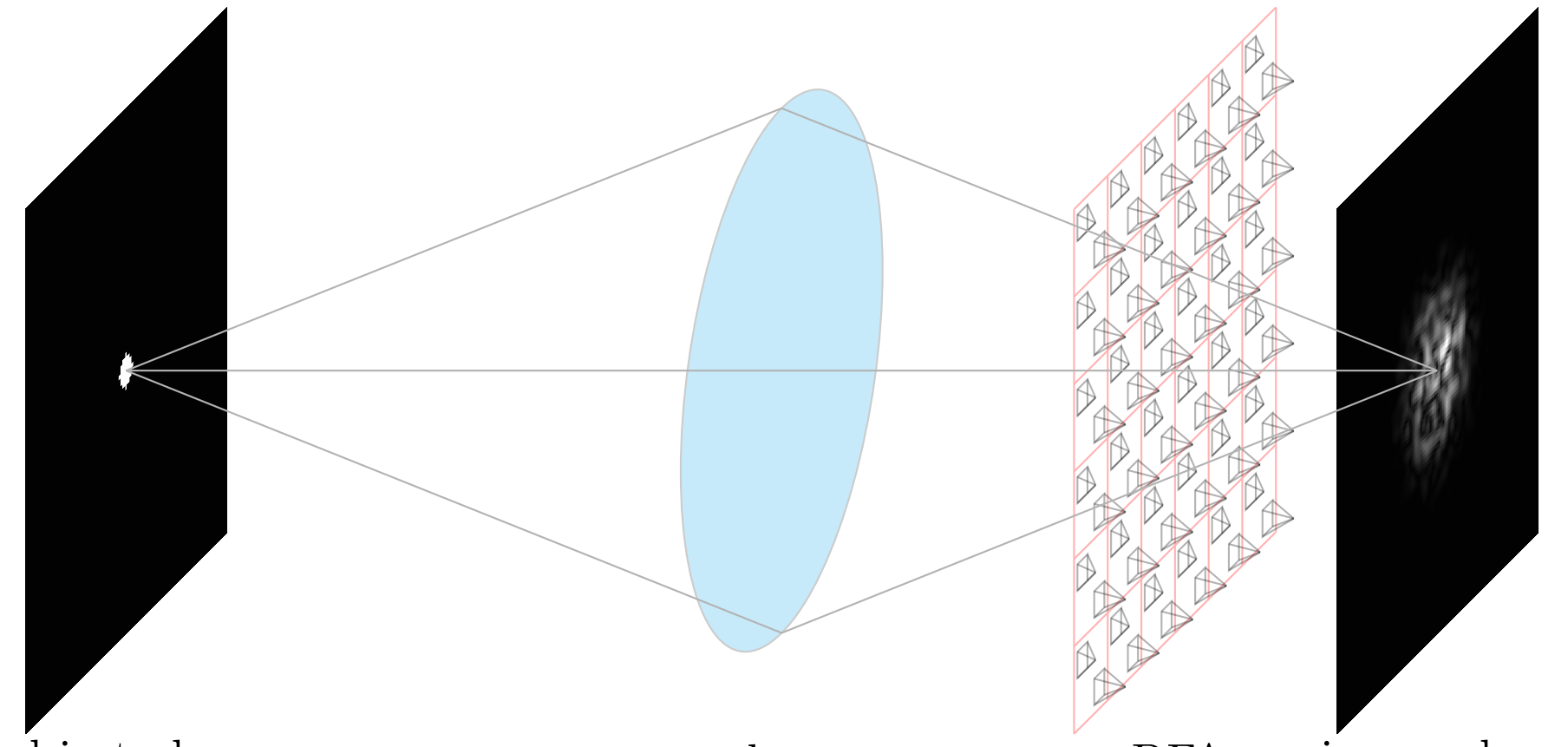
\begin{figure}
    \centering

    \begin{tikzpicture}[x=1cm,y=1cm]
      \begin{scope}[canvas is yz plane at x=-5]
       \node[transform shape] at (0,0)
         {\includegraphics[width=4cm,angle=-90]{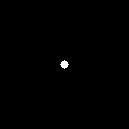}};
      \end{scope}
      \node at (-4.2,-2.2,2) {object plane};
      
      \begin{scope}[canvas is yz plane at x=0]
       \fill[cyan!20] (0,0) circle (2cm);
       \draw[black!20] (0,0) circle (2cm);
      \end{scope}
      \node at (0.8,-2.2,2) {lens};

      \begin{scope}[canvas is yz plane at x=5]
       \node[transform shape] at (0,0)
         {\includegraphics[width=4cm,angle=-90]{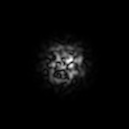}};
      \end{scope}
      \node at (5.8,-2.2,2) {image plane};
      
      \begin{scope}[canvas is yz plane at x=3]
       \foreach \i in {0,...,6} {
        \draw[red!30] (-2,-2+\i*4/6) -- (2,-2+\i*4/6);
       }
       \foreach \i in {0,...,6} {
        \draw[red!30] (-2+\i*4/6,-2) -- (-2+\i*4/6,2);
       }
      \end{scope}
      \begin{scope}[opacity=0.3]
      \foreach \i in {0,...,5} {
        \foreach \j in {0,...,5} {
	 \draw (3,-2+0.6*4/6+\i*4/6,-2+0.6*4/6+\j*4/6) -- (3,-2+0.9*4/6+\i*4/6,-2+0.6*4/6+\j*4/6) -- 
	       (3,-2+0.9*4/6+\i*4/6,-2+0.9*4/6+\j*4/6) -- (3,-2+0.6*4/6+\i*4/6,-2+0.9*4/6+\j*4/6) -- cycle;
	 \draw (3,-2+0.6*4/6+\i*4/6,-2+0.6*4/6+\j*4/6) -- (3.1,-2+0.75*4/6+\i*4/6,-2+0.75*4/6+\j*4/6);
	 \draw (3,-2+0.9*4/6+\i*4/6,-2+0.6*4/6+\j*4/6) -- (3.1,-2+0.75*4/6+\i*4/6,-2+0.75*4/6+\j*4/6);
	 \draw (3,-2+0.6*4/6+\i*4/6,-2+0.9*4/6+\j*4/6) -- (3.1,-2+0.75*4/6+\i*4/6,-2+0.75*4/6+\j*4/6);
	 \draw (3,-2+0.9*4/6+\i*4/6,-2+0.9*4/6+\j*4/6) -- (3.1,-2+0.75*4/6+\i*4/6,-2+0.75*4/6+\j*4/6);
	}
       }

      \foreach \i in {0,...,5} {
        \foreach \j in {0,...,5} {
	 \draw (3,-2+0.1*4/6+\i*4/6,-2+0.1*4/6+\j*4/6) -- (3,-2+0.4*4/6+\i*4/6,-2+0.1*4/6+\j*4/6) -- 
	       (3,-2+0.4*4/6+\i*4/6,-2+0.4*4/6+\j*4/6) -- (3,-2+0.1*4/6+\i*4/6,-2+0.4*4/6+\j*4/6) -- cycle;
	 \draw (3,-2+0.1*4/6+\i*4/6,-2+0.1*4/6+\j*4/6) -- (3.2,-2+0.25*4/6+\i*4/6,-2+0.25*4/6+\j*4/6);
	 \draw (3,-2+0.4*4/6+\i*4/6,-2+0.1*4/6+\j*4/6) -- (3.2,-2+0.25*4/6+\i*4/6,-2+0.25*4/6+\j*4/6);
	 \draw (3,-2+0.1*4/6+\i*4/6,-2+0.4*4/6+\j*4/6) -- (3.2,-2+0.25*4/6+\i*4/6,-2+0.25*4/6+\j*4/6);
	 \draw (3,-2+0.4*4/6+\i*4/6,-2+0.4*4/6+\j*4/6) -- (3.2,-2+0.25*4/6+\i*4/6,-2+0.25*4/6+\j*4/6);
	}
      }
      \end{scope}
      \node at (3.8,-2.2,2) {DFA};

      \draw[black!30] (-5,0,0) -- (5,0,0);
      \draw[black!30] (-5,0,0) -- (0,2,0) -- (5,0,0);
      \draw[black!30] (-5,0,0) -- (0,-2,0) -- (5,0,0);

    \end{tikzpicture}
 
    \caption{An optical system consisting of a lens, a periodic diffractive
    filter array (DFA) and a sensor array placed on the image plane. The DFA
    is made of a transparent material (say glass) with peaks that repeat
    periodically. The lens is assumed to be perfectly focused, i.e. the image of a
    point at the object plane is a point on the image plane. For illustration
    purposes, a point-like object (pinhole) is placed at the object plane and a
    cartoon point spread function at the image plane.}
    \label{fig:cartoon}
\end{figure}

The \emph{periodic circular convolution} or PCC naturally arises in optical
systems with a periodic element such as a periodic mask or diffractive
grating. In such systems, the point spread function or PSF (image of a point
through the optical system) depends on space in a predictable, periodic
fashion.  Convolution (in the usual sense) with the PSF is therefore 
inadequate to model such systems. 
An example of a system with a periodic PSF is given in
\cref{fig:cartoon}, where the PSF periodicity results from placing a
periodic diffractive filter array (DFA) close to the sensor. A similar
system is considered experimentally in \cite{Majumder:2025:HDS} with a
periodic DFA designed for hyperspectral imaging, i.e. the PSFs
are engineered to vary drastically with wavelength. The PCC algorithms
we present make it possible to solve the forward and inverse problem for such
systems efficiently, even if the PSFs spread over tens of pixels. The
application of PCCs to hyperspectral imaging is further explored in the
companion paper \cite{Gonzalez:2026:HID}. Here we solve a simpler,
single wavelength deconvolution problem with the PCC. The deconvolution
can be done efficiently because a matrix representation of the PCC is
unitarily equivalent to a block diagonal matrix. Hence common linear
algebra operations such as matrix multiplications and inversion are
efficient for PCCs.  The PCC naturally applies to
periodic difference operators on a finite lattice with periodic boundary
conditions, see e.g.  \cite{Berkolaiko:2013:IQG}. We use the block
diagonal representation of the PCC to compute a dispersion relation for
a periodic chain of masses, highlighting its relation to the
Floquet-Bloch transform.

\subsection{The PCC in one dimension}
\begin{figure}
\begin{center}
 \begin{tabular}{|c@{}c|c@{}c|}
  \hline
  $x$ & $K *_p x$ & $x$ & $K *_p x$\\\hline
  \includegraphics[width=0.2\textwidth]{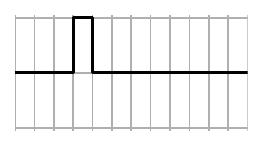} &
  \includegraphics[width=0.2\textwidth]{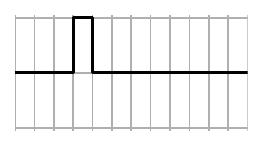} &
  \includegraphics[width=0.2\textwidth]{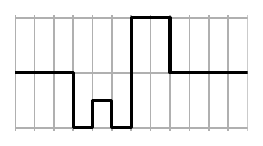} &
  \includegraphics[width=0.2\textwidth]{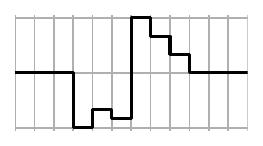} \\
  \includegraphics[width=0.2\textwidth]{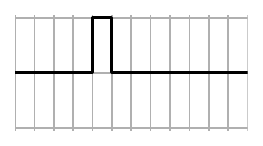} &
  \includegraphics[width=0.2\textwidth]{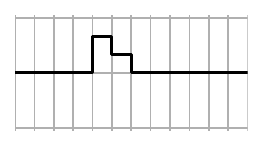} &
  \includegraphics[width=0.2\textwidth]{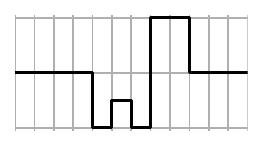} &
  \includegraphics[width=0.2\textwidth]{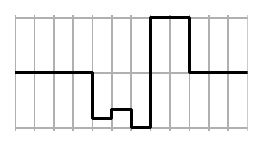} \\
  \includegraphics[width=0.2\textwidth]{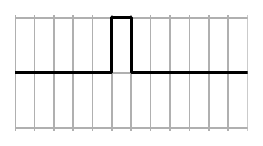} &
  \includegraphics[width=0.2\textwidth]{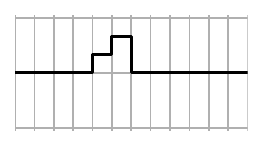} &
  \includegraphics[width=0.2\textwidth]{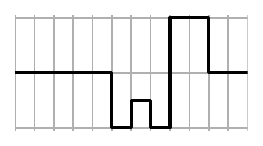} &
  \includegraphics[width=0.2\textwidth]{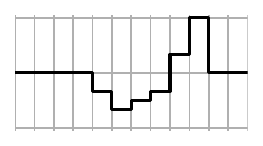}\\ \hline
 \end{tabular}
\end{center}

 \caption{Example of a one dimensional periodic circular
 convolution (PCC) with $n=12$ and $p=3$. The left pane visualizes the
 kernel $K$ by its effect on the canonical basis vectors $e_3$, $e_4$,
 and $e_5$. The PCCs of $e_0$, $e_1$, $e_2$ (resp. $e_6$, $e_7$,
 $e_8$) are similar but with $p=3$ circular shifts to the left (resp.
 right). The right pane shows the PCC of a signal and two
 circularly shifted versions.  As opposed to the circular convolution,
 the PCC of a circular shift is not the circular shift
 of the PCC. The gridlines for the vertical axis
 correspond to $-1$,$0$, and $1$.}
 \label{fig:pcc1d}
\end{figure}

A kernel $K$ acting on one dimensional signals in $\complex^n$ is an $n \times n$ second
order tensor with complex entries. Applying $K$ to a signal $x \in
\complex^n$ produces the signal $y \in \complex^n$ given by $y[i] =
\sum_{j=0}^{n-1} K[i,j] x[j]$, where $i=0,\ldots,n-1$ and tensor indices
start with zero. Thus kernel application is identical to
matrix vector multiplication for one dimensional signals. We say that $K$ is a \emph{one-dimensional
kernel with period $p$} if $p$ is a positive integer with $\mod{n}{p}=0$
and for any (signed) integers $i,j,\ell$  we have
\begin{equation}
 K[i + \ell p, j + \ell p ] = K[i,j],
 \label{eq:periodic:1d}
\end{equation}
where hereinafter tensor indices are understood
modulo the tensor size, e.g. in this case
$K[i,j]=K[\mod{i}{n},\mod{j}{n}]$ for any $i,j$ signed integers.
In particular, \eqref{eq:periodic:1d} means that column $j+p$ of $K$ is
column $j$ circularly shifted by $p$.
An example of a periodic one-dimensional kernel is given in
\cref{fig:pcc1d}.  Without exploiting the periodicity of $K$, applying
the kernel to an element of $\complex^n$ (a matrix vector product with
$K$) requires storing on the order of $n^2$ complex numbers and takes
$n^2$ floating point operations (flops). If the kernel has period $p$,
only $np$ numbers need to be stored.
We present algorithms that allow the kernel to be applied in
$\cO(n(p+\log n))$ flops, excluding setup costs. The special case of
period $p=1$ corresponds to the well-known circular convolution
\cite{Kincaid:1996:NA}. Indeed for any integers $i,j$, choose $\ell =
-j$ to see that 
\begin{equation}
 K[i,j] = K[i-j,0],
\end{equation}
so the kernel can be stored using the length $n$ vector $K[:,0]$ and can be
applied in $\cO(n\log n)$ flops. For general $p$, note that with the
Euclidean division of a signed integer $j$ by $p$, we can write $j$
uniquely as
$j=j'+j'' p$, where $j'=\mod{j}{p} \in \{0,1,\ldots,p-1\}$ and $j''$ is an integer.
Therefore for a kernel with period $p$ we have
\begin{equation}
 K[i,j] = K[i,j'+j''p] = K[i-j''p,j'],
\end{equation}
for any signed integers $i,j$. Hence, only the first $p$ columns of
$K$ are needed to apply $K$: the $n \times p$ matrix $K[:,0 \dv p-1]$ is
all that needs to be stored. We denote by $K *_p
x$ the application of a $p-$periodic kernel $K$ to a signal $x$, with the
understanding that $*_1 \equiv *$ is the circular convolution operation. The PCC can be
computed using the naïve 
\cref{alg:naive1}, which costs as much as a matrix
vector product ($n^2$ flops) but reduces the storage to $np$
entries of $K$. For the sake of conciseness we use the
following notation. For $i\leq j$, we denote an integer range as
$\irange{i,j} = \{ i,i+1,\ldots,j\}$, in that order.  Similarly, an
index progression with a constant positive step $s$ is denoted by $i \dv
s \dv j = \{i, i+s, \ldots, i+ms\}$, where $m$ is the largest integer
such that $i+ms \leq j$. If $s$ is omitted, then $s=1$ is assumed, i.e.
$i \dv j = \irange{i,j}$. A colon in an index denotes the whole range of
that index, i.e. $K[:,0]$ is the first column of $K$ and $\hat K[:,i'',:]$
is the $i''-$th slice of $\hat K$. The vector of all zeroes of size $n$ is
denoted by $0_n$ and the matrix of all zeroes of size $n_1 \times
n_2$ is $0_{n_1 \times n_2}$. Similarly $1_n$ and $1_{n_1 \times n_2}$
denote a vector or matrix of all ones.

\begin{algorithm}[h]
\caption{Naïve 1D PCC.}
\label{alg:naive1}
\textbf{Inputs:} A signal $x \in \complex^n$ and $K[:,0 \dv p-1] \in
\complex^{n\times p}$ for a $p-$periodic kernel $K$.

\textbf{Outputs:} A signal $y \in \complex^n$ s.t. $y = K *_p x$.

\begin{tabbing}
xxx\= xxx \= \kill
$y=0_n$\\
\textbf{for} $i \in \irange{0,n-1}$, $j \in \irange{0,n-1}$\\
\> $y[i] = y[i] + x[j] K[i-j + (\mod{j}{p}), \mod{j}{p}]$\\
\textbf{end}\\
\textbf{return} $y$
\end{tabbing}
\end{algorithm}

\subsection{The PCC in two dimensions}
\begin{figure}
 \begin{tabular}{|@{\hspace{0.10em}}c@{\hspace{0.10em}}|@{\hspace{0.10em}}c@{\hspace{0.10em}}c@{\hspace{0.10em}}c@{\hspace{0.10em}}c@{\hspace{0.10em}}c@{\hspace{0.10em}}c@{\hspace{0.10em}}c@{\hspace{0.10em}}c@{\hspace{0.10em}}c|}
 
  \hline
  \raisebox{1.5em}{\rotatebox{90}{$x$}} &
  \rule{0pt}{0.11\textwidth}%
  \includegraphics[width=0.10\textwidth]{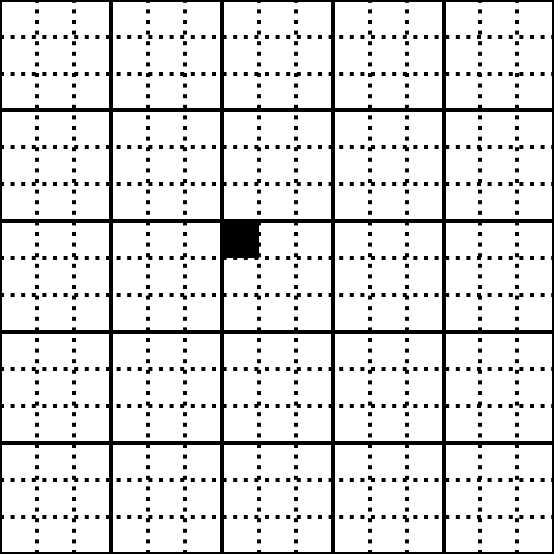}&
  \includegraphics[width=0.10\textwidth]{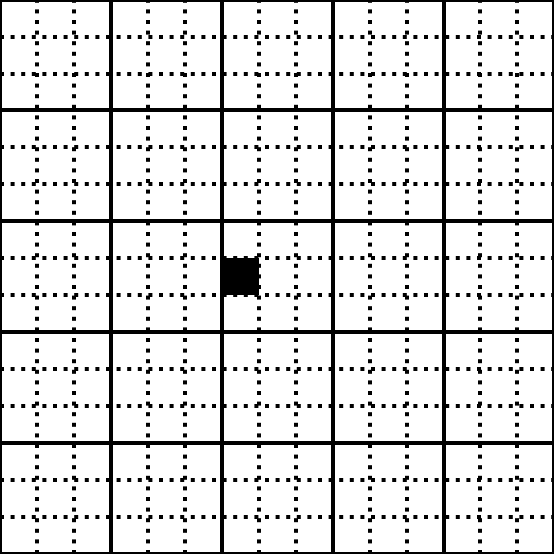}& 
  \includegraphics[width=0.10\textwidth]{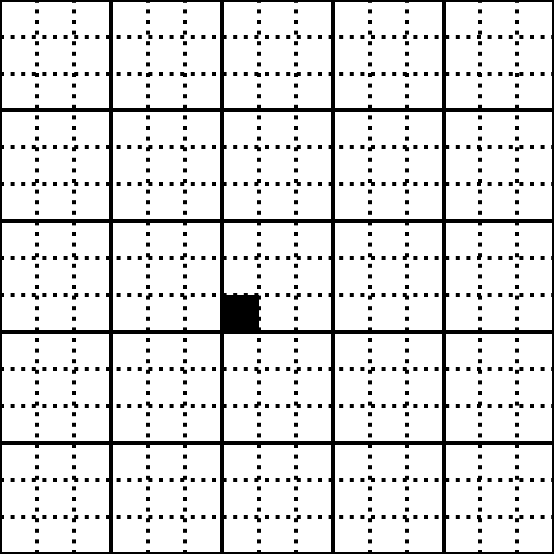}& 
  \includegraphics[width=0.10\textwidth]{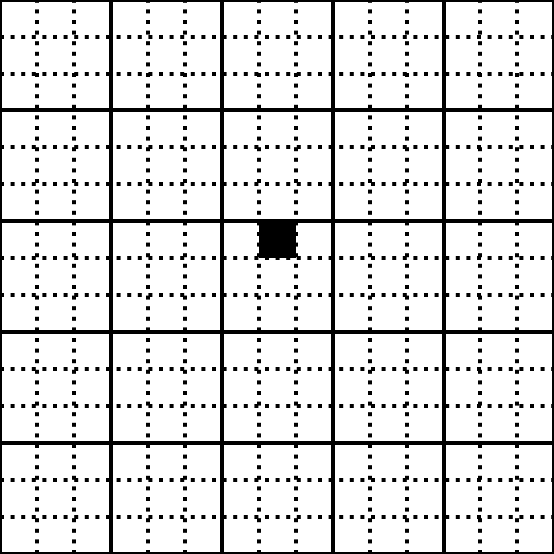}& 
  \includegraphics[width=0.10\textwidth]{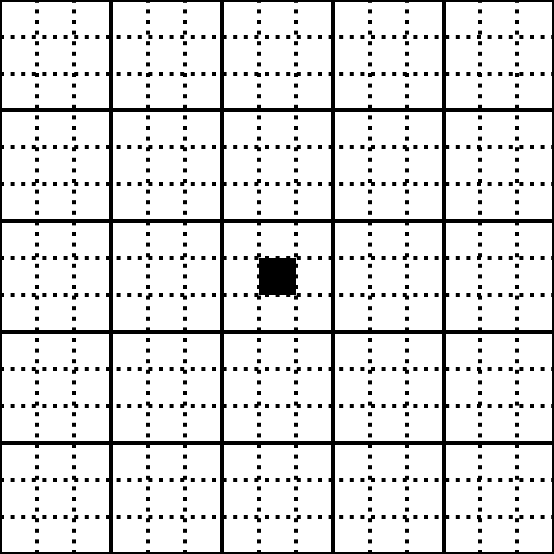}& 
  \includegraphics[width=0.10\textwidth]{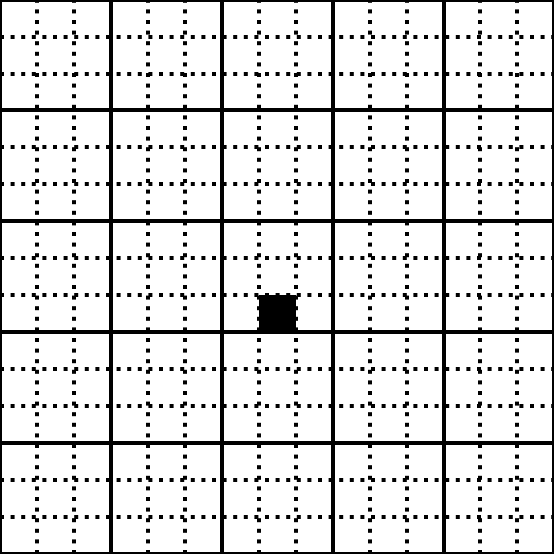}& 
  \includegraphics[width=0.10\textwidth]{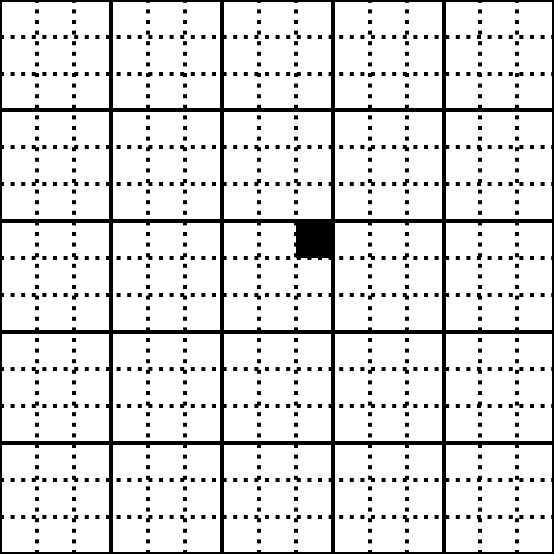}& 
  \includegraphics[width=0.10\textwidth]{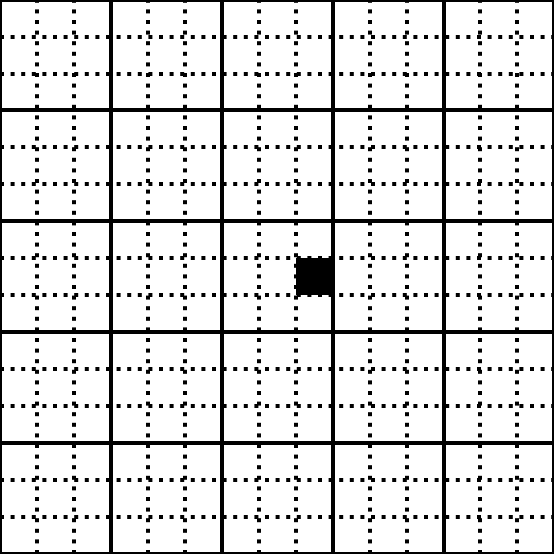}& 
  \includegraphics[width=0.10\textwidth]{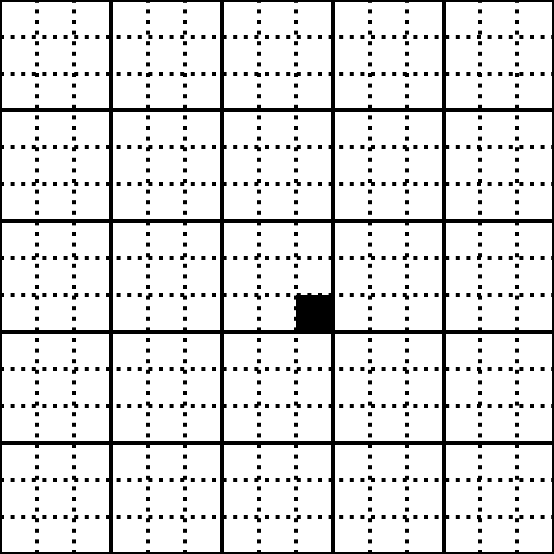} \\[-0.3em]
  \raisebox{1em}{\rotatebox{90}{$K *_p x$}} &
  \includegraphics[width=0.10\textwidth]{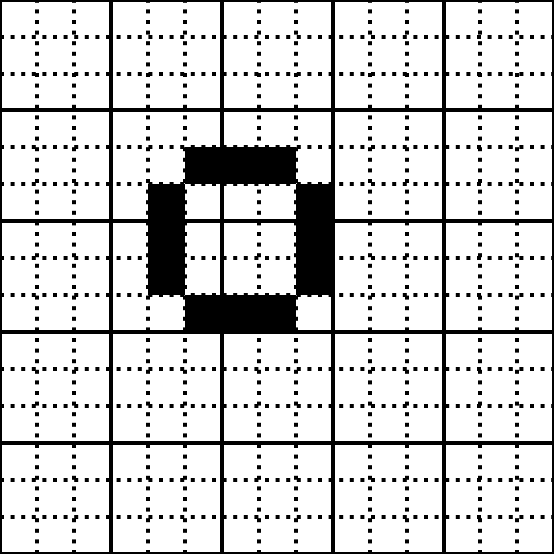}& 
  \includegraphics[width=0.10\textwidth]{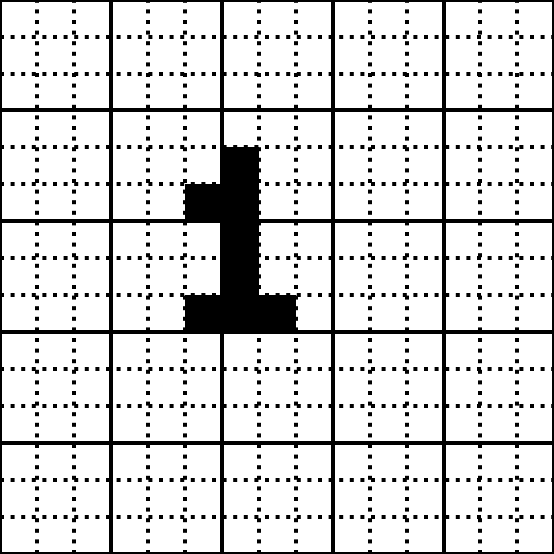}& 
  \includegraphics[width=0.10\textwidth]{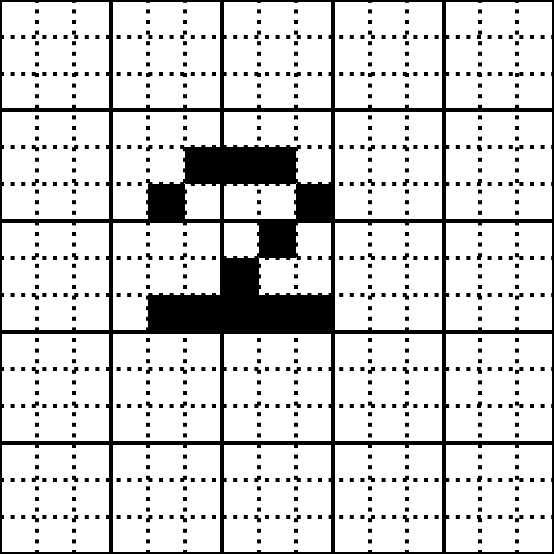}& 
  \includegraphics[width=0.10\textwidth]{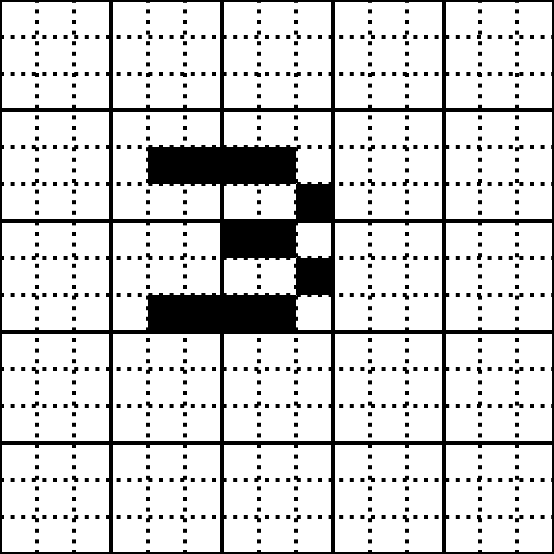}& 
  \includegraphics[width=0.10\textwidth]{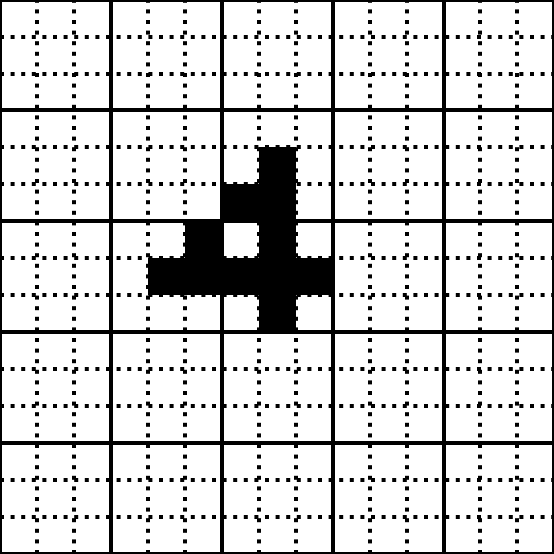}& 
  \includegraphics[width=0.10\textwidth]{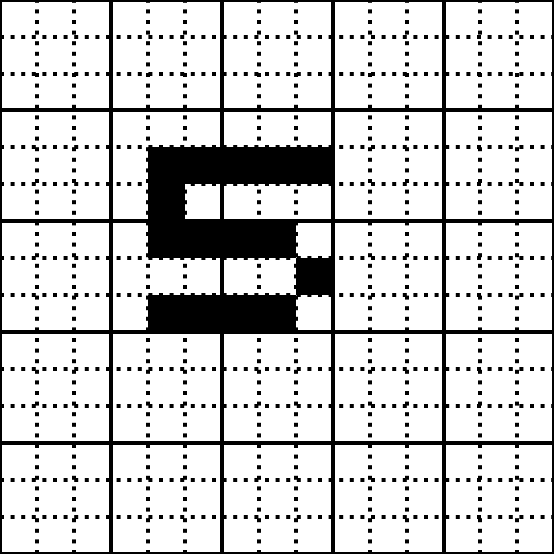}& 
  \includegraphics[width=0.10\textwidth]{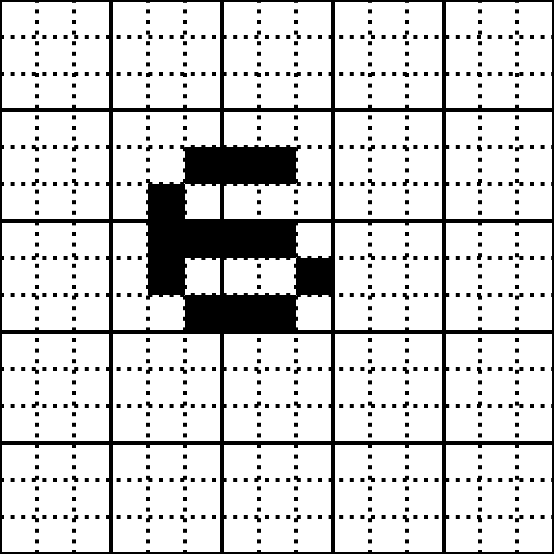}& 
  \includegraphics[width=0.10\textwidth]{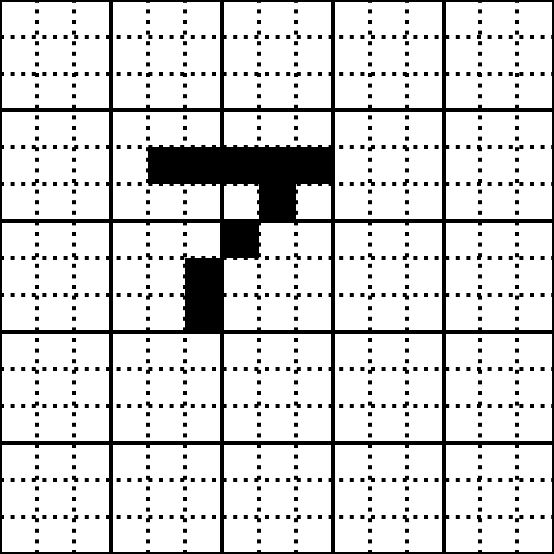}& 
  \includegraphics[width=0.10\textwidth]{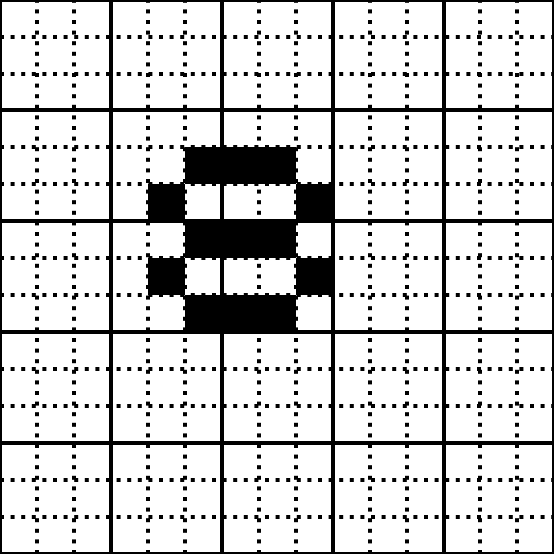} \\\hline
  \raisebox{1.5em}{\rotatebox{90}{$x$}} &
  \rule{0pt}{0.11\textwidth}%
  \includegraphics[width=0.10\textwidth]{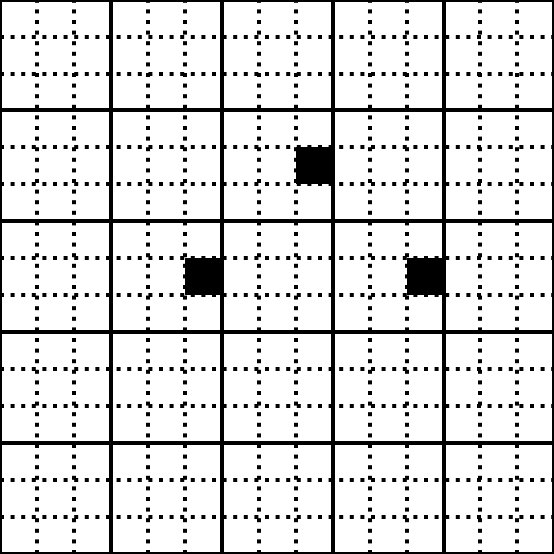}& 
  \includegraphics[width=0.10\textwidth]{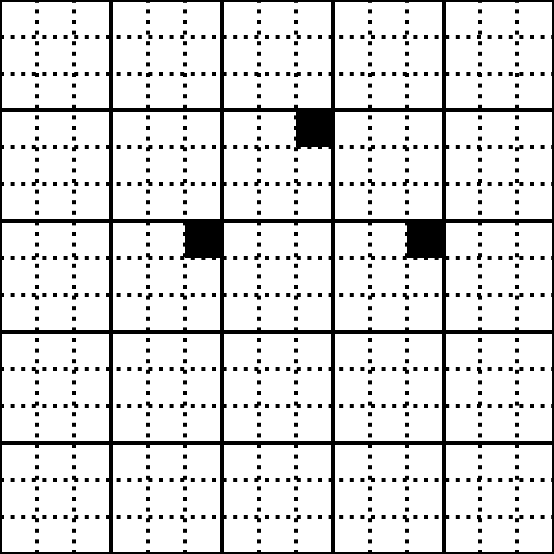}& 
  \includegraphics[width=0.10\textwidth]{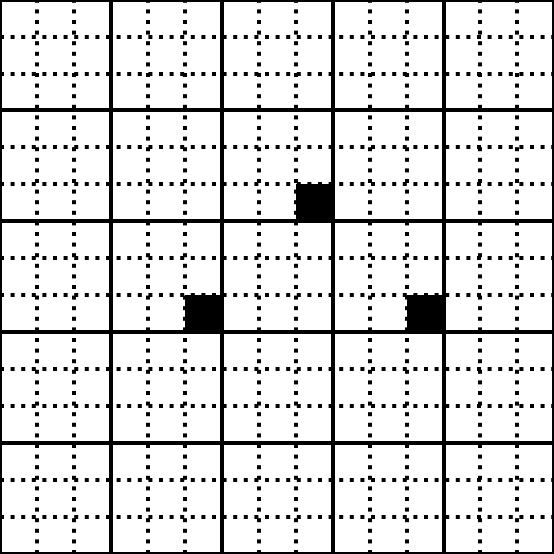}& 
  \includegraphics[width=0.10\textwidth]{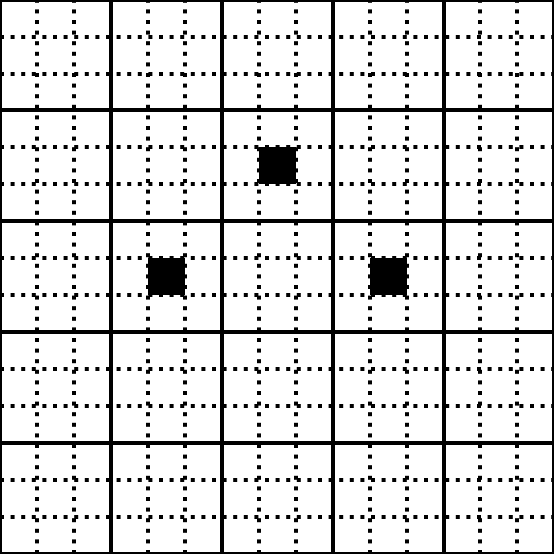}& 
  \includegraphics[width=0.10\textwidth]{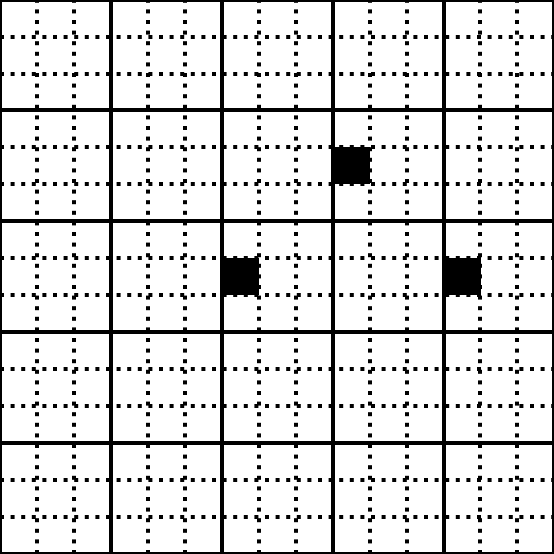}& 
  \includegraphics[width=0.10\textwidth]{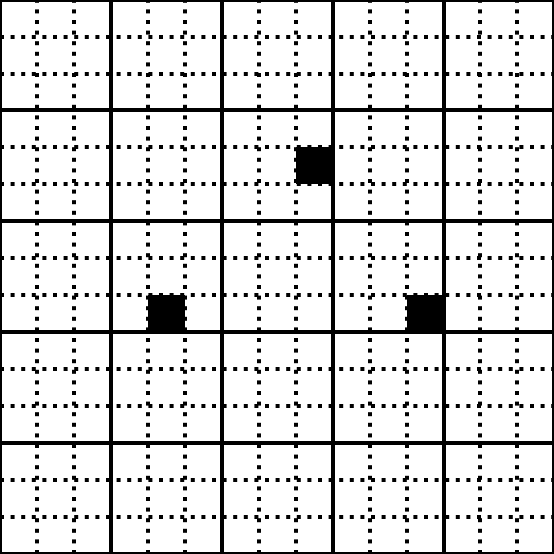}& 
  \includegraphics[width=0.10\textwidth]{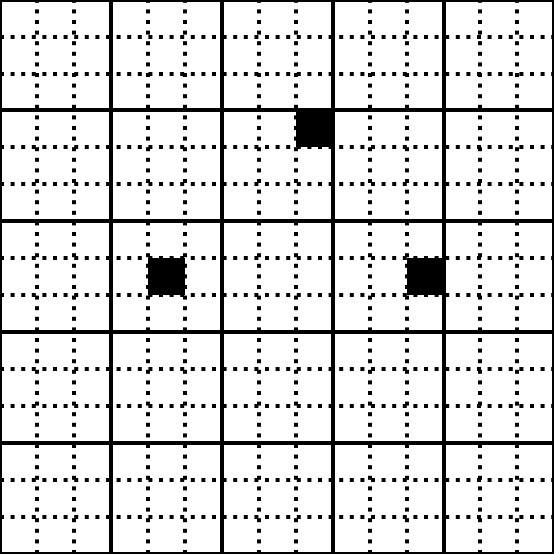}& 
  \includegraphics[width=0.10\textwidth]{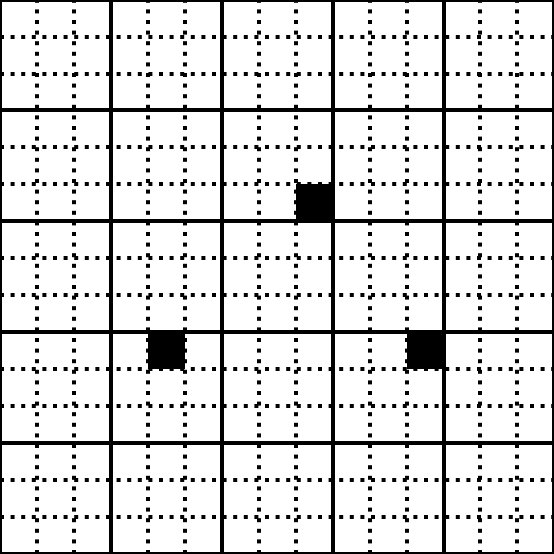}& 
  \includegraphics[width=0.10\textwidth]{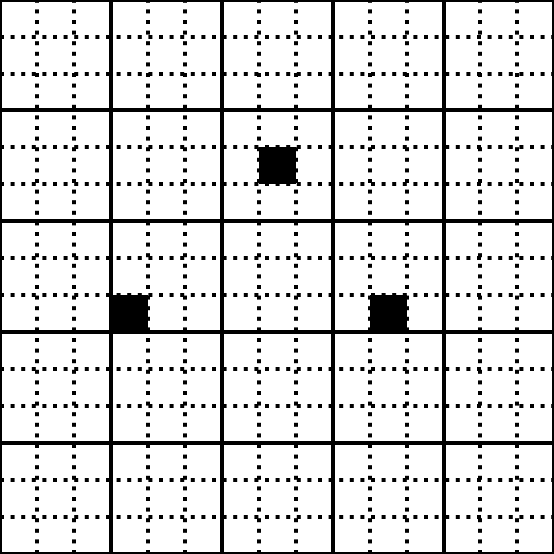} \\[-0.3em]
  \raisebox{1em}{\rotatebox{90}{$K *_p x$}} &
  \includegraphics[width=0.10\textwidth]{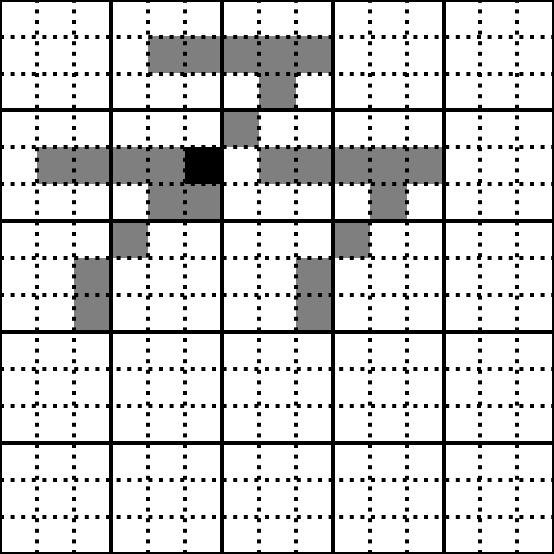}& 
  \includegraphics[width=0.10\textwidth]{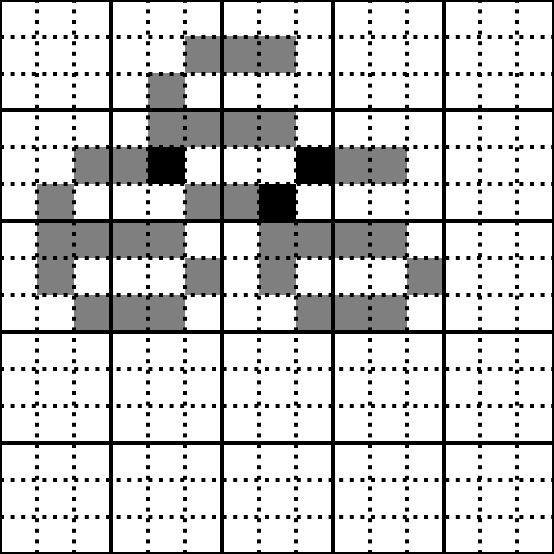}& 
  \includegraphics[width=0.10\textwidth]{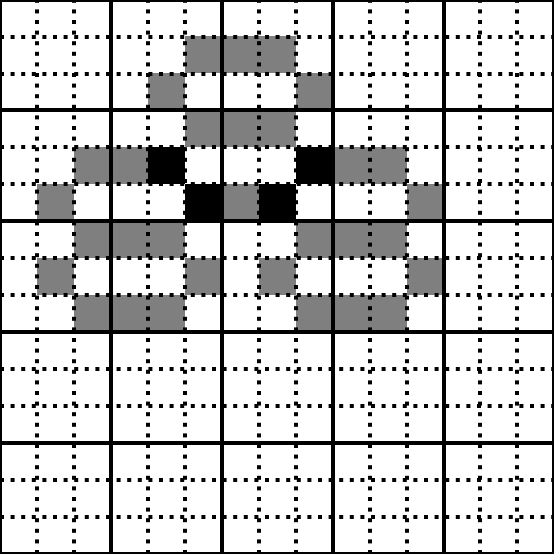}& 
  \includegraphics[width=0.10\textwidth]{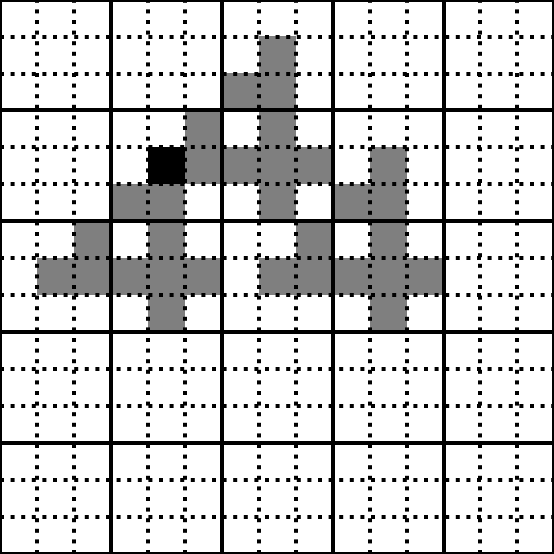}& 
  \includegraphics[width=0.10\textwidth]{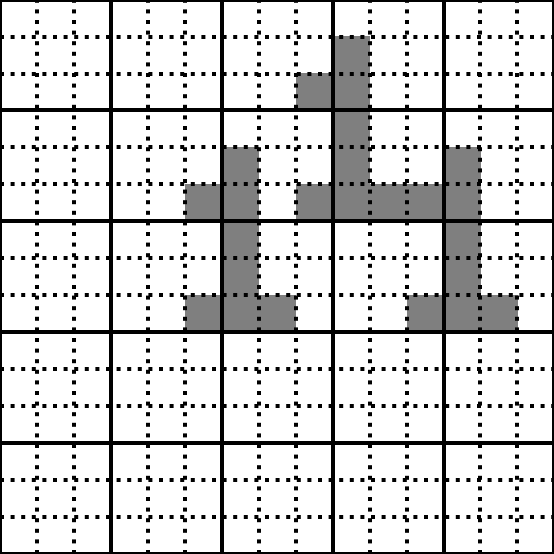}& 
  \includegraphics[width=0.10\textwidth]{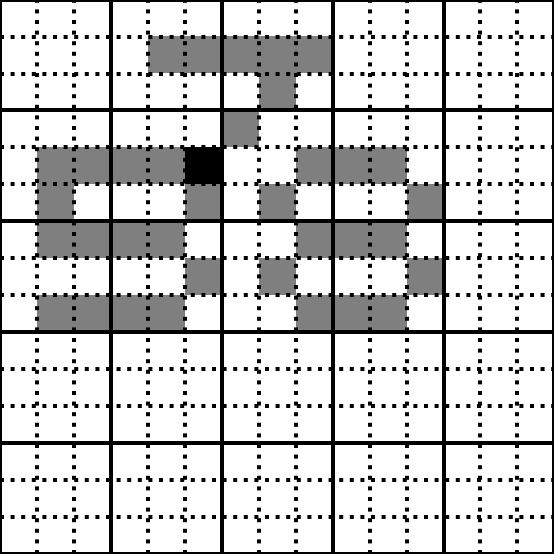}& 
  \includegraphics[width=0.10\textwidth]{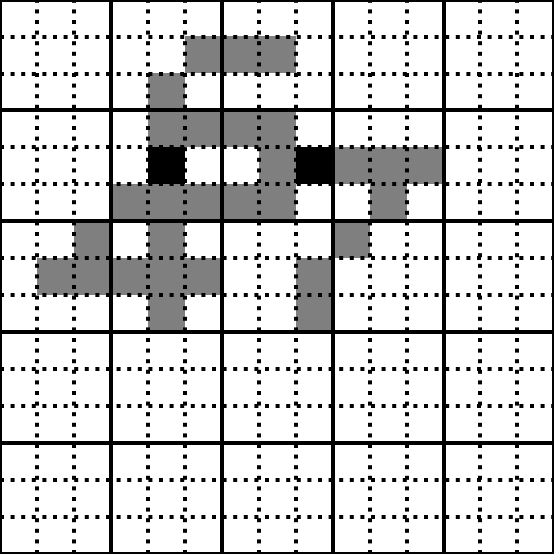}& 
  \includegraphics[width=0.10\textwidth]{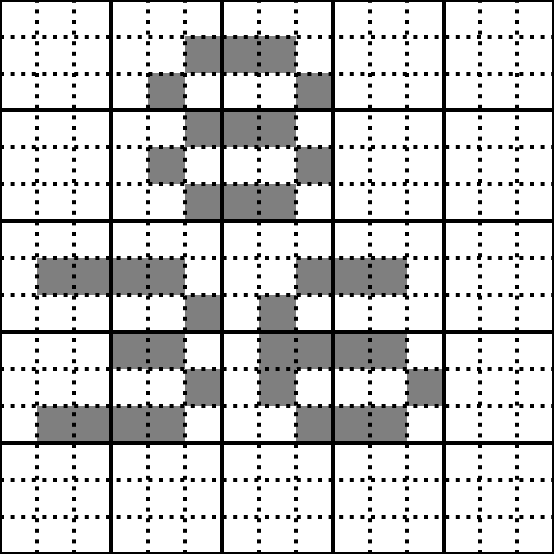}& 
  \includegraphics[width=0.10\textwidth]{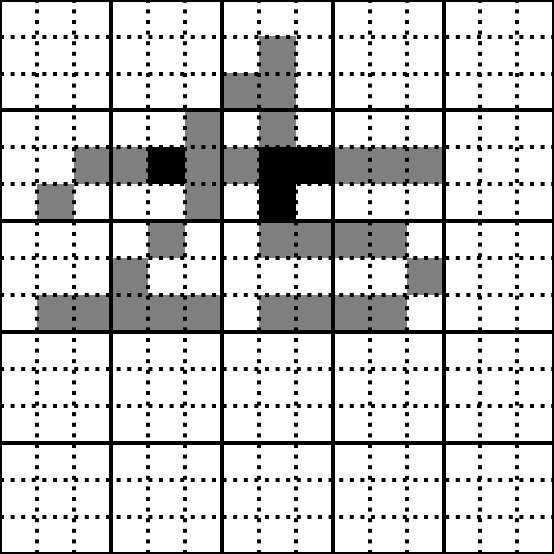} \\[0em]\hline
 \end{tabular}
 \caption{Illustration of 2D PCC with $n=(15,15)$ and $p=(3,3)$. The
 kernel $K$ is illustrated in the top pane by its PCC on canonical basis
 vectors. The kernel corresponds to digits from $0$ to $8$ to be
 recognizable by the human eye and is not meant to model any physical
 system. The bottom pane shows the PCC of different signals with the
 same kernel. The first five signals in the bottom pane are composed of
 three pixels separated by integer multiples of $p$ and thus display
 three copies of the same digit. The last four signals are also composed
 of three pixels, but since the spacings are not congruent with the
 period, they show three different digits. The solid lines represent
 periods. The color scale is the same on each row and ranges between
 white (zero) and black (maximum).}
 \label{fig:pcc2d}
\end{figure}

A kernel $K$ acting on two dimensional signals in $\complex^{n_1 \times
n_2}$ is a $4-$th order tensor of dimensions $n_1 \times n_2 \times n_1
\times n_2$. Applying $K$ to a signal $x \in \complex^{n_1 \times n_2}$
produces the signal $y \in \complex^{n_1 \times n_2}$ given by
\begin{equation}
 y[i_1,i_2] = \sum_{j_1=0}^{n_1-1} \sum_{j_2=0}^{n_2-1}
 K[i_1,i_2,j_1,j_2] x[j_1,j_2],
\end{equation}
for any $i_1,i_2$ signed integers.
Letting $n=(n_1,n_2)$ and $|n|=n_1n_2$, we see that applying a
two-dimensional kernel requires $|n|^2$ operations, while storing
$|n|^2$ scalars for the kernel. We say that $K$ is a
\emph{two-dimensional kernel with 
period $p=(p_1,p_2)$} if $p_i$ is a positive integer such that $\mod{n_i}{p_i} = 0$,
$i=1,2$ and if for any integers $i_1,i_2,j_1,j_2,\ell_1,\ell_2$ we have
\begin{equation}
  K[i_1 + \ell_1 p_1, i_2 + \ell_2 p_2, j_1 + \ell_1 p_1, j_2 + \ell_2
  p_2] = K[i_1,i_2,j_1,j_2].
 \label{eq:periodic:2d}
\end{equation}
As in the 1D case, we observe that only $|p||n|$ entries need to be
stored for a kernel $K$ with period $p$. Indeed for any
$i_1,i_2,j_1,j_2$ integers,
\begin{equation}
 K[i_1,i_2,j_1,j_2] = K[i_1-j_1''p_1, i_2-j_2''p_2, j_1', j_2' ].
\end{equation}
where we used Euclidean division to write $j_m = j_m' + j_m'' p_m$,
with $j_m' \in \irange{0,p_m-1}$, $j_m'' \in \irange{0,q_m-1}$, $n_m =
p_m q_m$ and $m=1,2$.
A naïve algorithm can be written which stores only
$|n||p|$ entries rather than $|n|^2$ but still takes $|n|^2$ operations
to compute. An example of the 2D PCC operation is given in \cref{fig:pcc2d}.
The case $|p|=1$ also corresponds to the usual
circular convolution of 2D signals. We present here fast algorithms that
can apply the kernel in $\cO(|p||n| + |n|\log|n|)$ flops. The extension
to $d-$dimensional signals is possible with  $2d-$th order tensors as
kernels, however this is omitted for brevity.

\subsection{Contents}
We start by recalling a few facts about the discrete Fourier transform
(DFT) in \cref{sec:not}. A first PCC algorithm that reduces the cost to
about $|p|+1$ discrete Fourier transforms (DFTs), excluding setup costs,
is presented in \cref{sec:pcc1}. Since the cost of a DFT is $\cO(|n|
\log |n|)$, this is already a significant improvement over the naïve
implementation. In \cref{sec:pcc2} we reduce the cost further to about
$2$ DFTs, independently of the period $p$, by showing that the PCC
operation is up to DFTs equivalent to a tensor contraction over $d$
indices where $d=1$ or $2$ is the dimension of the problem.  The
computational cost in terms of floating point operations is estimated in
\cref{sec:cost}. A numerical comparison of both approaches is presented
in \cref{sec:num}. Crucially, the tensor contraction can be
represented, up to permutations, by a block diagonal matrix.  Hence PCCs
can be operated and analyzed with the same ease as block matrices, e.g.
they can be multiplied and inverted efficiently (\cref{sec:blk}).  We
present two example applications to deconvolution \cref{sec:deconv} and
to dispersion relation calculations \cref{sec:greenper}.  Summary and
conclusions appear in \cref{sec:conc}.

\section{The discrete Fourier transform}
\label{sec:not}
Our results rely on the discrete Fourier transform or DFT. Recall that 
the DFT of a 1D signal $x \in \complex^n$
is the signal $\cF_n x \in \complex^n$ defined by (see e.g.
\cite{Kincaid:1996:NA})
\begin{equation}
 (\cF_n x)[i] = \sum_{j=0}^{n-1} \omega_n^{ij} x[j],~i=0,\ldots,n-1,
 \label{eq:dft}
\end{equation}
where $\omega_n = \exp[-\imath 2\pi / n ]$ and $\imath = \sqrt{-1}$. 
The circular convolution of signals in $\complex^n$ can be computed
efficiently in $\cO(n \log n)$ operations using the DFT:
\begin{equation}
 \label{eq:conv}
 \cF_n (K * x) = ( \cF_n x ) \hprod ( \cF_n K[:,0] ),
\end{equation}
where $\hprod$ is the componentwise or Hadamard product. 
The DFT of a second order tensor $x \in \complex^{n_1 \times n_2}$ is given by
\begin{equation}
 (\cF_{n_1,n_2} x )[i_1,i_2] = \sum_{j_1=0}^{n_1-1} \sum_{j_2=0}^{n_2-1}
  \omega_{n_1}^{i_1j_1} \omega_{n_2}^{i_2j_2}x[j_1,j_2].
\end{equation} 
The circular convolution of two tensors of the same size $n_1 \times
n_2$ can be computed as in \eqref{eq:conv}, replacing $\cF_n$ by
$\cF_{n_1,n_2}$, at cost of $\cO(|n|\log |n|)$ flops. The DFT and
circular convolution can be defined for higher order tensors in a
similar fashion. Finally, if one of the dimensions in the DFT is $*$, we
mean that there is no DFT along that direction. For example for $x \in
\complex^{n_1 \times n_2}$,
\begin{equation}
 (\cF_{*,n_2} x)[i_1,i_2] =\sum_{j_2=0}^{n_2-1}
   \omega_{n_2}^{i_2j_2}x[i_1,j_2].
\end{equation}

\section{PCC as circular convolutions of size $n$}
\label{sec:pcc1}
We first show that a 1D PCC with a $p-$periodic kernel can be computed
using $p$ circular convolutions. To see this, let us first consider the
case of signals $x\in \complex^n$ that have only non-zero components at
the indices $\ell \dv p \dv n-1$ (i.e.  indices $i$ with $\mod{i}{p} =
\ell$). For such signals, sometimes called \emph{Dirac combs}, the PCC
reduces to a circular convolution. 
\begin{lemma}
 \label{lem:comb1}
 Let $K \in \complex^{n\times n}$ be a $p$ periodic kernel. Let
 $\ell\in\irange{0,p-1}$ and consider $x \in \complex^n$ such that
 all the entries of $x$ are zero except possibly at the indices
 $\ell \dv p \dv n-1$. Then we have
 \[
 K *_p x = K[:+\ell,\ell] * x,
 \]
where  $K[:+\ell,\ell]$ is a circular shift of $K[:,\ell]$ to the
\emph{left} by $\ell$.
\end{lemma}
\begin{proof}
Let $q$ be such that $n=pq$. By using the zero pattern of $x$ we see
that
\begin{equation}
 (K *_p x)[i] = \sum_{j=0}^{n-1} K[i,j] x[j] = \sum_{j=0}^{q-1}
 K[i,\ell+pj] x[\ell+pj]
\end{equation}
Since $K$ is $p-$periodic we have
\begin{equation}
 (K *_p x)[i] = \sum_{j=0}^{q-1} K[i-pj,\ell] x[\ell+pj] =
 \sum_{j=0}^{q-1} K[i-(\ell+pj)+\ell,\ell] x[\ell+pj].
\end{equation}
Using the zero pattern of $x$ again, we recognize a circular convolution
to conclude:
\begin{equation}
 (K *_p x)[i] = (K[:+\ell,\ell]*x)[i].
\end{equation}
\end{proof}

An immediate consequence is that a one-dimensional PCC with a
$p-$periodic kernel $K$ can be calculated using $p$ circular
convolutions. 
\begin{theorem}
\label{thm:pcc1}
 For $x \in \complex^n$ we have
 \begin{equation}
  K *_p x = \sum_{i=0}^{p-1} K[:+i,i] * x_i,
 \end{equation}
 where $x_i[i \dv p \dv n-1]=x[i \dv p \dv n-1]$ and $x_i$ is otherwise zero.
\end{theorem}
\begin{proof}
Since $x=x_0 + \cdots + x_{p-1}$, we may use linearity of the PCC and
\cref{lem:comb1} to obtain the desired result.
\end{proof}

By using the convolution theorem $p$ times we obtain \cref{alg:dft1} for
computing the PCC of 1D signals. Clearly, the result of the preparation step can
be saved if the same kernel $K$ is applied to many
signals $x$. A sample  Julia \cite{Bezanson:2017:JFA} implementation is
included in \texttt{PCC.jl} \cite{PCC} with the pair of functions
\verb|pcc1d_gemv_cc_prep| and \verb|pcc1d_gemv_cc|.

\begin{algorithm}
\caption{PCC with size $n$ circular convolutions.}
\label{alg:dft1}
\textbf{Inputs:} A signal $x \in \complex^n$ and $K[:,0\dv p-1] \in
\complex^{n\times p}$ for a $p-$periodic kernel $K$.

\textbf{Outputs:} A signal $y \in \complex^n$ s.t. $y = K *_p x$.

\begin{tabbing}
 xxx\= xxx \= \kill
 \texttt{\# Preparation}\\
 $\hat K = 0_{n \times p}$\\
 \textbf{for} $i\in \irange{0,p-1}$\\
 \> $\hat K[:,i] = \cF_n K[:+i,i]$\\
 \textbf{end}\\
 \texttt{\# Application}\\
 $y=0_n$\\
 \textbf{for} $i\in \irange{0,p-1}$\\
 \> $x_i= 0_n$ \\
 \> $x_i[i \dv p \dv n-1] = x[i \dv p \dv n-1]$\\
 \> $y = y + \hat K [:,i] \hprod (\cF_n x_i)$\\
\textbf{end}\\
\textbf{return} $\cF_n^{-1} y$
\end{tabbing}
\end{algorithm}

A similar result can be proven for 2D signals and is included below. The proof
is omitted as it is very similar to that of \cref{thm:pcc1}.
\begin{theorem}
\label{thm:pcc2}
 For some $x \in \complex^{n_1 \times n_2}$ the PCC with a
 $p=(p_1,p_2)-$periodic kernel $K \in \complex^{n_1 \times n_2 \times
 n_1 \times n_2}$ can be calculated as follows
 \begin{equation}
  K *_p x = \sum_{i_1=0}^{p_1-1} \sum_{i_2=0}^{p_2-1}
  K[:+i_1,:+i_2,i_1,i_2] * x_{i_1,i_2},
 \end{equation}
 where $x_{i_1,i_2}[i_1 \dv p_1 \dv n_1 -1 , i_2 \dv p_2 \dv n_2 - 1] = x[i_1
 \dv p_1 \dv n_1 -1, i_2 \dv p_2 \dv n_2 - 1]$ and $x_{i_1,i_2}$ is otherwise
 zero.
\end{theorem}

\section{PCC as tensor contraction}
\label{sec:pcc2}
The algorithms in \cref{sec:pcc1} rely on expressing the PCC as $|p|$
circular convolutions of size $n$. We show here that the PCC can
also be written as $|p|^2$ circular convolutions of size $q$, where $q$
is the integer such that $n=pq$ for 1D signals and $q=(q_1,q_2)$ where
$n_i=p_iq_i$, $i=1,2$ for 2D signals. This alternate interpretation of the PCC
is advantageous as it requires less operations and shows that the PCC is up to
DFTs a tensor contraction over $d$ indices, where $d$ is the signal's dimension.
We present the derivation for the 1D PCC in detail in
\cref{sec:ipcc1} and then present in \cref{sec:ipcc2} the 2D result,
omitting some details for brevity. 

\subsection{One dimensional PCC}
\label{sec:ipcc1}
We start by showing that the 1D PCC with a $p-$periodic kernel can
be written as the sum of $p^2$ circular convolutions of size $q$.
\begin{lemma}
 \label{lem:pcc1:q}
 Let $K$ be an $n\times n$ one-dimensional kernel with period $p$ and let
 $q$ be such that $n=pq$. Then for any $x \in \complex^n$ we have
 \begin{equation}
  \label{eq:pcc1:q}
  K *_p x = \sum_{i'=0}^{p-1} \sum_{j'=0}^{p-1} 
   K[i'\dv p \dv n-1, j'] * x[j'\dv p \dv n-1],
 \end{equation}
 where each of the circular convolutions is of size $q$.
\end{lemma}
\begin{proof}
Let $i=i'+p i'' \in \irange{0,n-1}$, where $i' \in \irange{0,p-1}$ and $i'' \in
\irange{0,q-1}$. Similarly we write $j \in \irange{0,n-1}$ as
$j=j'+pj''$. With these new indices, the application of the kernel $K$ to $x$ becomes
\begin{equation}
 (K *_p x)[i'+pi''] = \sum_{j'=0}^{p-1} \sum_{j''=0}^{q-1} K[i'+p i'',
 j' + p j''] x[j'+ p j''].
\end{equation}
But $K$ is $p-$periodic so we can write
\begin{equation}
 (K *_p x)[i'+pi''] =  \sum_{j'=0}^{p-1} \sum_{j''=0}^{q-1} K[i'+p
 (i''-j''), j'] x[j'+ p j''].
\end{equation}
For each $i',j'$ fixed, the sum over $j''$ is a circular convolution of size $q$.
Therefore
\begin{equation}
 (K *_p x)[i' \dv p \dv n-1] = \sum_{j'=0}^{p-1} K[i' \dv p \dv n-1,j'] * x[j' \dv p \dv
 n-1].
\end{equation}
The result follows by summing over $i'$ to obtain all the entries of $K
*_p x$.
\end{proof}

Inspecting the right hand side of \cref{eq:pcc1:q}, we see that the
convolutions of size $q$ inside the sum are accessing entries of $K$ and
$x$ with stride $p$. This may be better explained if we make an
assumption on the memory ordering of entries in a tensor. We assume here that
all tensors are stored in \emph{column-major ordering}, as in Matlab,
Julia or Fortran, meaning that when enumerating the entries of a tensor,
the first (leftmost) index is the index that changes the fastest. We
find it convenient to use ``reshape'' operations to express the access
pattern in \cref{eq:pcc1:q}.  For example, $\tilde K = \reshape (K[:,0
\dv p-1],p,q,p)$ takes the entries in the $n \times p$ second order
tensor $K[:,0 \dv p-1]$ (which are stored column by column), and
reinterprets them as a $p \times q \times p$ third order tensor $\tilde
K$.

With this in mind, we use  the convolution theorem for each of the size
$q$ circular convolutions in \cref{lem:pcc1:q}, to obtain the following
generalization of the convolution theorem \cref{eq:conv} to PCCs.  The
following result shows that a PCC is, up to an appropriate pair of DFTs,
a tensor contraction operation.  
\begin{theorem}
\label{thm:gconv1}
Let $K$ be an $n \times n$, $p-$periodic kernel, $x \in
\complex^n$ and $y = K *_p x$. Define the reshaped tensors
\begin{equation}
\begin{aligned}
\tilde K &= \reshape(K[:,0 \dv p-1],p,q,p),\\
\tilde x &= \reshape(x,p,q)~\text{and}\\
\tilde y &= \reshape(y,p,q)
\end{aligned} 
\end{equation}
with $q$ such that $n=pq$. Then for $i'\in \irange{0,p-1}$ and $i''\in
\irange{0,q-1}$ we have
\begin{equation}
 \label{eq:gconv1}
 \hat y [i',i''] = \sum_{j'=0}^{p-1} \hat K[i',i'',j'] \hat x[j',i''],
\end{equation}
where $\hat K = \cF_{*,q,*} \tilde K$, $\hat x = \cF_{*,q} \tilde x$
and $\hat y = \cF_{*,q} \tilde y$.
\end{theorem}
\begin{proof}
First note that for any $v \in \complex^n$, we have
$\reshape(v,p,q)[i,j] = v[i+pj]$, where $i \in \irange{0,p-1}$ and $j\in
\irange{0,q-1}$. We can rewrite
\eqref{eq:pcc1:q} as follows 
\begin{equation}
 \tilde y[i',:] = \sum_{j'=0}^{p-1} \tilde K[i',:,j'] * \tilde x[j',:],
 ~\text{for}~i' \in \irange{0,p-1}.
\end{equation}
By the convolution theorem \eqref{eq:conv}, we get
\begin{equation}
 \hat y[i',:] = \sum_{j'=0}^{p-1} \hat K[i',:,j'] \hprod \hat x[j',:],
\end{equation}
which written in index notation gives the desired result
\eqref{eq:gconv1}.
\end{proof}
\begin{remark}
There are two extreme cases for $p$.
First when $p=1$, \cref{thm:gconv1} reduces to the usual convolution theorem
\eqref{eq:conv}. Second when $p=n$, all DFTs become the identity and
\eqref{eq:gconv1} reduces to a matrix vector product with same computational cost
as the naïve algorithm \cref{alg:naive1}.
\end{remark}

An implementation of the PCC using \cref{thm:gconv1} is presented in
\cref{alg:idft1}.  Similarly to \cref{alg:dft1}, the calculations carried
out during the preparation step can be stored and used many times, if
the same kernel is applied to calculate the PCC of many different
signals. A sample implementation is included in \verb|PCC.jl| with the
pair of functions \verb|pcc1d_gemv_improved_prep| and
\verb|pcc1d_gemv_improved|.  

\begin{algorithm}
 \caption{Improved DFT PCC algorithm in 1D.}
\label{alg:idft1}

\textbf{Inputs:} A signal $x \in \complex^n$ and $K[:,0\dv p-1] \in
\complex^{n\times p}$ for a $p-$periodic kernel $K$.

\textbf{Outputs:} A signal $y \in \complex^n$ s.t. $y = K *_p x$.
\begin{tabbing}
xxx\= xxx \= \kill
\texttt{\# Preparation}\\
$\hat K = \cF_{*,q,*}\reshape(K[:,0\dv p-1],p,q,p)$\\
\texttt{\# Application}\\
$y = 0_{p \times q}$\\
$\hat x = \cF_{*,q} \reshape(x,p,q)$\\
\textbf{for} $\ell \in \irange{0,p-1}$, $j\in\irange{0,q-1}$, $i \in \irange{0,p-1}$\\
\> $y[i,j] = y[i,j] + \hat K[i,j,\ell] \hat x [\ell,j]$\\
\textbf{end}\\
\textbf{return} $\reshape(\cF_{*,q}^{-1} y,n)$
\end{tabbing}

\end{algorithm}

\subsection{Two dimensional PCC}
\label{sec:ipcc2}
We start by showing that the PCC of a kernel with period
$p=(p_1,p_2)$ and a signal in $\complex^{n_1 \times n_2}$ is the sum of
$|p|^2$ circular convolutions of size $q=(q_1,q_2)$, where
$n_i=p_iq_i$, $q_i$ integer for $i=1,2$. The proof of
\cref{lem:pcc2:q} is omitted because of its similarity to that of
\cref{lem:pcc1:q}.
\begin{lemma}
 \label{lem:pcc2:q}
 Let $K$ be an $n_1 \times n_2 \times n_1 \times n_2$ 2D kernel with
 period $p=(p_1,p_2)$ and let
 $q=(q_1,q_2)$ be such that $n_i=p_iq_i$, $i=1,2$. Then for any $x \in
 \complex^{n_1\times n_2}$ we have
 \begin{equation}
  \label{eq:pcc2:q}
  \begin{aligned}
  K *_p x = \sum_{i_1'=0}^{p_1-1}& \sum_{j_1'=0}^{p_1-1} 
\sum_{i_2'=0}^{p_2-1} \sum_{j_2'=0}^{p_2-1}\\
   &K[i_1'\dv p_1 \dv n_1-1, i_2'\dv p_2 \dv n_2-1, j_1',j_2'] *
   x[j_1'\dv p_1 \dv n_1-1,j_2'\dv p_2 \dv n_2-1],
 \end{aligned}
\end{equation}
 where each of the circular convolutions is of size $q$.
\end{lemma}
As in the 1D case, tensor notation can be used to apply the
convolution theorem $|p|^2$ times in order to evaluate \eqref{eq:pcc2:q}
efficiently. The proof of the next \cref{thm:gconv2} is also omitted
because of its similarity to the proof of \cref{thm:gconv1}.
\begin{theorem}
\label{thm:gconv2}
 Let $K$ be an $n_1 \times n_2 \times n_1 \times n_2$ 2D kernel with period $p=(p_1,p_2)$.
 Let $x \in \complex^{n_1 \times n_2}$ and $y = K *_p x$.  Define 
 \begin{equation}
 \begin{aligned}
  \tilde K &=\reshape(K[:,:,0 \dv p_1-1,0 \dv
 p_2-1],p_1,q_1,p_2,q_2,p_1,p_2), \\
  \tilde x &= \reshape(x,p_1,q_1,p_2,q_2)~\text{and}\\
 \tilde y &= \reshape(y,p_1,q_1,p_2,q_2). 
 \end{aligned}
 \end{equation}
Then for $i_l'\in \irange{0,p_l-1}$ and $j_l\in
\irange{0,q_l-1}$, $l=1,2$, we have
\begin{equation}
 \label{eq:gconv2}
 \hat y [i_1',i_1'',i_2',i_2''] = \sum_{j_1'=0}^{p_1-1}
 \sum_{j_2'=0}^{p_2-1}
 \hat K[i_1',i_1'',i_2',i_2'',j_1',j_2'] \hat x[j_1',i_1'',j_2',i_2''],
\end{equation}
where $\hat K = \cF_{*,q_1,*,q_2,*,*} \tilde K$, $\hat x =
\cF_{*,q_1,*,q_2} \tilde x$ and  $\hat y = \cF_{*,q_1,*,q_2} \tilde y$.
\end{theorem}

The improved PCC algorithm in 2D is included in \cref{alg:idft2}. The
preparation phase involves computing $|p|^2$ DFTs of size $q$. The
application phase involves computing $2|p|$ DFTs of size $q$ and a
tensor contraction on two indices of a $6-$th order tensor with a $4-$th
order one.   An inverse Fourier transform finalizes the computation. A
sample implementation in Julia is included in \texttt{PCC.jl} as the
pair of functions \verb|pcc2d_gemv_improved_prep| and
\verb|pcc2d_gemv_improved|, where the nested for loop is different
than in \cref{alg:idft2} and ordered so that the traversal of $\hat K$
(the largest object in memory) is linear in memory.

\begin{algorithm}
\caption{Improved DFT PCC algorithm in 2D.}
\label{alg:idft2}
\textbf{Inputs:} A signal $x \in \complex^{n_1\times n_2}$ and a
$K[:,:,0 \dv p_1-1, 0 \dv p_2-1] \in \complex^{n_1 \times n_2 \times p_1 \times
p_2}$, for a $p-$periodic kernel $K$.

\textbf{Outputs:} A signal $y \in \complex^{n_1 \times n_2}$ s.t. $y = K
*_p x$.
\begin{tabbing}
xxx\= xxx \= \kill
\texttt{\# Preparation}\\
$\hat K = \cF_{*,q_1,*,q_2,*,*}
\reshape(K[:,:,0 \dv p_1-1,0 \dv p_2-1],p_1,q_1,p_2,q_2,p_1,p_2)$\\
\texttt{\# Application}\\
$y = 0_{p_1,q_1,p_2,q_2}$\\
$\hat x = \cF_{*,q_1,*,q_2} \reshape(x,p_1,q_1,p_2,q_2)$\\
\textbf{for} $i_1 \in \irange{0,p_1-1},~
	        j_1 \in \irange{0,q_1-1},~
                i_2 \in \irange{0,p_2-1},~
	        j_2 \in \irange{0,q_2-1}$\\
\>\textbf{for} $\ell_1 \in \irange{0,p_1-1},~
                \ell_2 \in \irange{0,p_2-1}$\\
\>\> $y[i_1,j_1,i_2,j_2] = y[i_1,j_1,i_2,j_2] + 
\hat K[i_1,j_1,i_2,j_2,\ell_1,\ell_2] \hat x [\ell_1,j_1,\ell_2,j_2]$\\
\>\textbf{end}\\
\textbf{end}\\
\textbf{return} $\reshape(\cF_{*,q_1,*,q_2}^{-1} y,n_1,n_2)$
\end{tabbing}
\end{algorithm}

\begin{remark}
\label[remark]{rem:blk}
In both 1D and 2D, the PCC can be represented by a matrix that is
unitarily equivalent to a block diagonal matrix with blocks of size $|p|
\times |p|$. This can be readily seen in 1D: the sum in
\eqref{eq:gconv1} calculates the
$p \times p$ matrix vector products $\hat K[:,i'',:] \hat x[:,i'']$ for $i'' \in
\irange{0,q-1}$. In 2D, the matrix vector products from \eqref{eq:gconv2} are
\begin{equation}
\reshape(\hat K[:,i_1'',:,i_2'',:,:],|p|,|p|) \reshape(\hat
x[:,i_1'',:,i_2''],|p|).
\end{equation}
This block matrix product approach is implemented in \texttt{PCC.jl} as
the functions \verb|pcc?d_gemv_blkm_prep| and \verb|pcc?d_gemv_blkm|,
where the question mark is the dimension, either 1 or 2. Their
computational cost is, ignoring memory access costs, identical to that
of \cref{alg:idft2}. In practice we observed wall clock timings similar
to those of \cref{alg:idft2}. The block diagonal matrix point of view
is crucial to carry common linear algebra operations on PCCs, see
\cref{sec:blk} for more details. An illustration for one dimensional
PCCs is included in \cref{fig:blkmat}.
\end{remark}

\setlength{\imagesize}{0.2\textwidth}
\newcommand{\amatrix}[1]{%
  \mathord{\vcenter{\hbox{%
    \tikz[x=\dimexpr\imagesize/5\relax,
          y=-\dimexpr\imagesize/5\relax]{
          \path[use as bounding box] (0,-0.1) rectangle (5,5.1);
	  \node[anchor=north west, inner sep=0] at (0,0)
            {\includegraphics[width=\imagesize]{#1}};
          \foreach \i in {1,2,3,4} {
           \draw[thick,blue] (0,\i)--(5,\i);
           \draw[thick,blue] (\i,0)--(\i,5);
          }
	  }%
  }}}%
}

\newcommand{\smatrix}[1]{%
  \mathord{\vcenter{\hbox{%
    \tikz[x=\dimexpr\imagesize/5\relax,
          y=-\dimexpr\imagesize/5\relax]{
          \path[use as bounding box] (0,-0.1) rectangle (1,5.1);
	  \node[anchor=north west, inner sep=0] at (0,0)
            {\includegraphics[height=\imagesize]{#1}};
          \foreach \i in {1,2,3,4} {
           \draw[thick,blue] (0,\i)--(1,\i);
          }
	  }%
  }}}%
}

\begin{figure}
 \begin{center}
 \begin{tikzpicture}[x=\dimexpr\imagesize/5\relax,
                     y=-\dimexpr\imagesize/5\relax]
 \node at (0,0)
 {$
  K = \left[ \amatrix{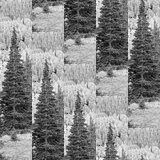} \right] 
  = 
  \cF_{*,q}^{-1} 
  \left[ \amatrix{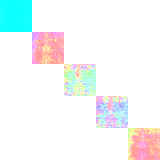} \right] 
  \cF_{*,q}
 $};
 
 \draw[->,thick] (-6,3) -- node[midway,above] {$\cF_{*,q}$} (-1,6);
 \draw[->,thick] (1,6) -- (4,3);
 \node at (0,6) {$ \left[ \smatrix{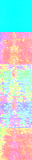} \right] $};
 \draw[red,thick, rounded corners=5pt,opacity=0.5] (-6.6,2.6) rectangle (-5.4,-2.6);
 \end{tikzpicture}
 \end{center}
 \caption{An illustration of \cref{rem:blk} for a one-dimensional PCC
 $K$ with period $p$ and dimension $n \times n$ with $n=pq$. Here $p=32$
 and $q=5$.  Blocks of dimension $32 \times 32$ are highlighted and the
 gray color scale goes from 0 (white) to 1 (black). The periodicity of
 $K$ is apparent because a block-column is a circular shift by $p=32$ of
 the previous block-column. Therefore we only need to store the entries
 in the first block-column (highlighted in red). The DFT $\cF_{*,q}$ of
 the first \mbox{block-column} is computed to obtain a block diagonal
 representation of the matrix $K$. Colors display the
 complex phases stemming from the DFT.}
 \label{fig:blkmat}
\end{figure}

\section{Computational cost}
\label{sec:cost}
The different PCC algorithms we have presented can be separated into
\emph{preparation} and \emph{application} phases. During the preparation
phase, we carry out computations whose results can be stored in memory
and reused many times during the application phase, to compute the PCC
of a fixed kernel with many different signals. We make a few simplifying
assumptions to estimate the computational cost. First, the cost of a DFT
of size $n$ is approximated by $|n| \log |n|$ flops, where the base of
the logarithm is fixed, say base $2$, and $n$ could be either an integer
for the 1D DFT or a tuple $(n_1,n_2)$ for the 2D DFT. This cost estimate
is based on the Fast Fourier Transform algorithm of Cooley and Tukey
\cite{Cooley:1965:AMC}, see also \cite{Johnson:2026:IFP}. Finally we
track only the number of multiplication operations and we ignore costs
associated with memory access and allocation. Memory costs are
important, as we see in the numerical results of \cref{sec:num}.

The naïve implementation \cref{alg:naive1} uses (part of) the kernel
as-is, so there is no preparation phase involved. To apply
the PCC with the naïve algorithm involves $|n|^2$ flops.
This is also true for the 2D version, although we do not include the
implementation in the main text. 

For the first DFT implementation \cref{alg:dft1} (or its two dimensional
variant included as \texttt{pcc2d\_fft} in \texttt{PCC.jl}), the
preparation phase involves computing the DFT of $|p|$ signals of size $n$. Hence we
estimate its cost to be $|p| |n| \log |n|$ flops. Applying the first DFT
implementation involves $|p|$ DFTs of signals with size $|n|$ followed
by $|p||n|$ multiplications. The result is
finalized by an additional DFT. Hence
the total cost is about $(|p|+1)|n|\log|n|$ flops.

The preparation phase for the improved DFT 
\cref{alg:idft1,alg:idft2} involves computing $|p|^2$ DFTs of size $q$,
which amounts to $|p|^2 |q|\log |q| = |p| |n| (\log |n| - \log|p|)$ flops, a slightly
lower cost than that of the first DFT implementation.
The application phase is significantly cheaper. It can be decomposed in
three steps. The first and third steps each consist of $|p|$
DFTs of size $q$, with a cost of $|p| |q| \log |q| = |n| (\log |n| -
\log|p|)$ flops. The second step is a tensor contraction involving
$|p|^2|q| = |p||n|$ flops. Hence the improved algorithm costs $|p||n| +
2 |n| (\log |n| - \log|p|)$ flops.

The computational cost estimates are summarized in \cref{tab:cost}. By
dividing the application cost of the improved vs convolution algorithms
we see that the expected speedup for $|p|$ fixed approaches $1 +
\log|n|$ for large $|n|$.

\begin{table}
\caption{Computational cost estimate in floating point operations
(flops) for each of the PCC algorithms. For a fixed kernel, the
preparation operations need to be computed only once and the result can
then be reused to compute the PCC with many different signals. The
application cost corresponds to the cost per signal. For reference the algorithms are the naïve (\cref{alg:naive1}), DFT (\cref{alg:dft1}) and DFT improved (\cref{alg:idft1,alg:idft2}).}
\label{tab:cost}
\begin{center}
\begin{tabular}{|l|l|l|}
\hline
 PCC  & Preparation & Application\\\hline
 Naïve  & 0 & $|n|^2$\\
 Convolution & $|p| |n| \log |n|$ & $|p||n| + (|p|+1) |n| \log |n|$\\
 Improved & $|p| |n| (\log |n|-\log|p|) $ & $|p||n| + 2 |n| (\log
 |n|-\log|p|)$\\\hline
\end{tabular}
\end{center}
\end{table}

\section{Benchmarks}
\label{sec:num}
\begin{figure}
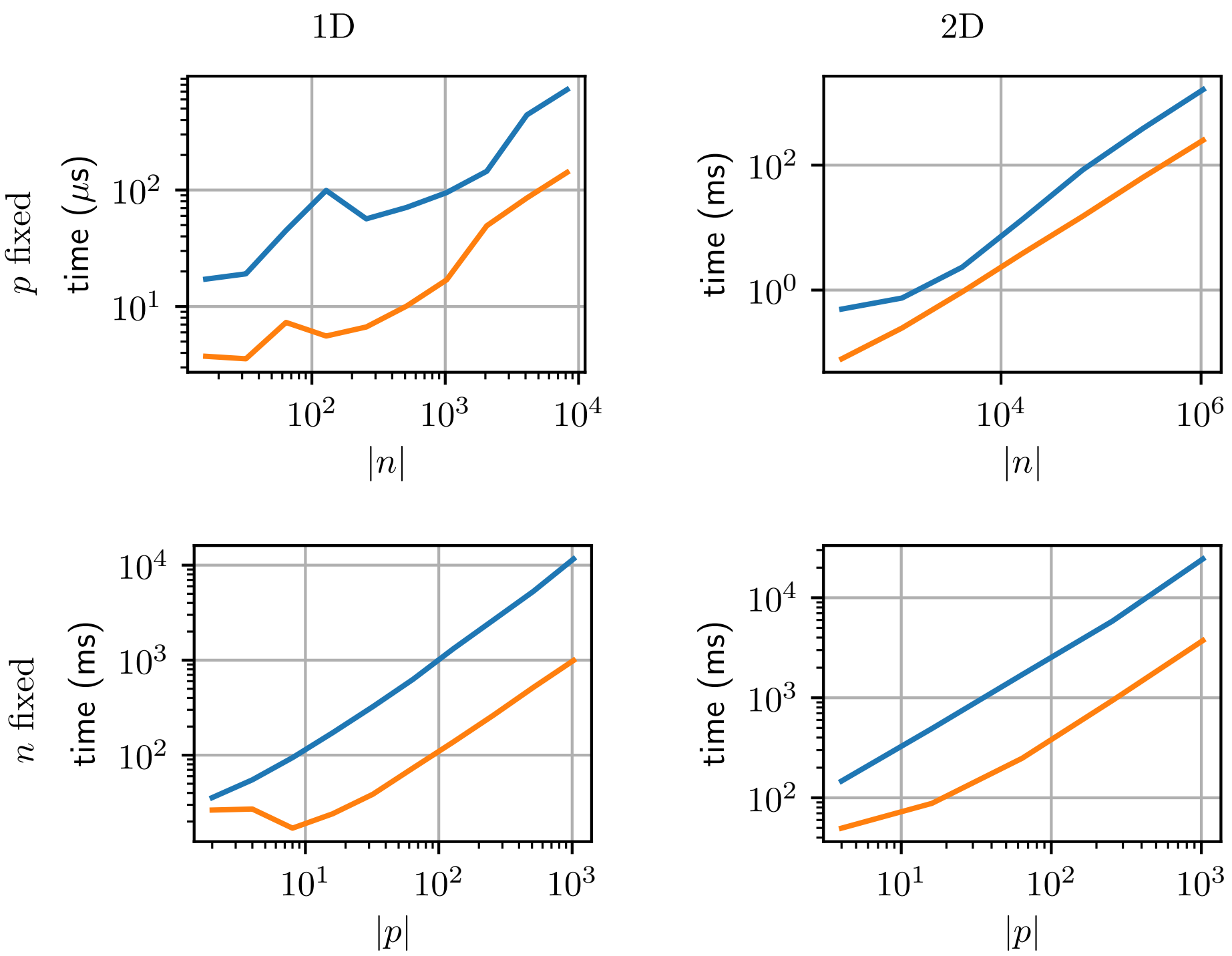

 \begin{tabular}{c@{}cc}
  & 1D & 2D\\
  \raisebox{6em}{\rotatebox{90}{$p$ fixed}} &
  \input{pfixed_time_1d.pgf} & 
  \input{pfixed_time_2d.pgf}\\
  \raisebox{6em}{\rotatebox{90}{$n$ fixed}} &
  \input{nfixed_time_1d.pgf} & 
  \input{nfixed_time_2d.pgf}
 \end{tabular}
 \caption{Wall-clock times for implementations of the PCC with circular
 convolutions (in orange) and as tensor contraction (the ``improved''
 algorithm, in blue).
 The first row corresponds to $p$ fixed experiments, where  $p=8$ for 1D
 or $p=(8,8)$ for 2D and $n$ varies.  For the $n$ fixed experiment we
 used $n=2^{10}$ in 1D and $n=(2^{10},2^{10})$ in 2D and varied $p$. For
 more details on the sizes considered see \cref{tab:summary}.}
 \label{fig:time}
\end{figure}

\begin{figure}
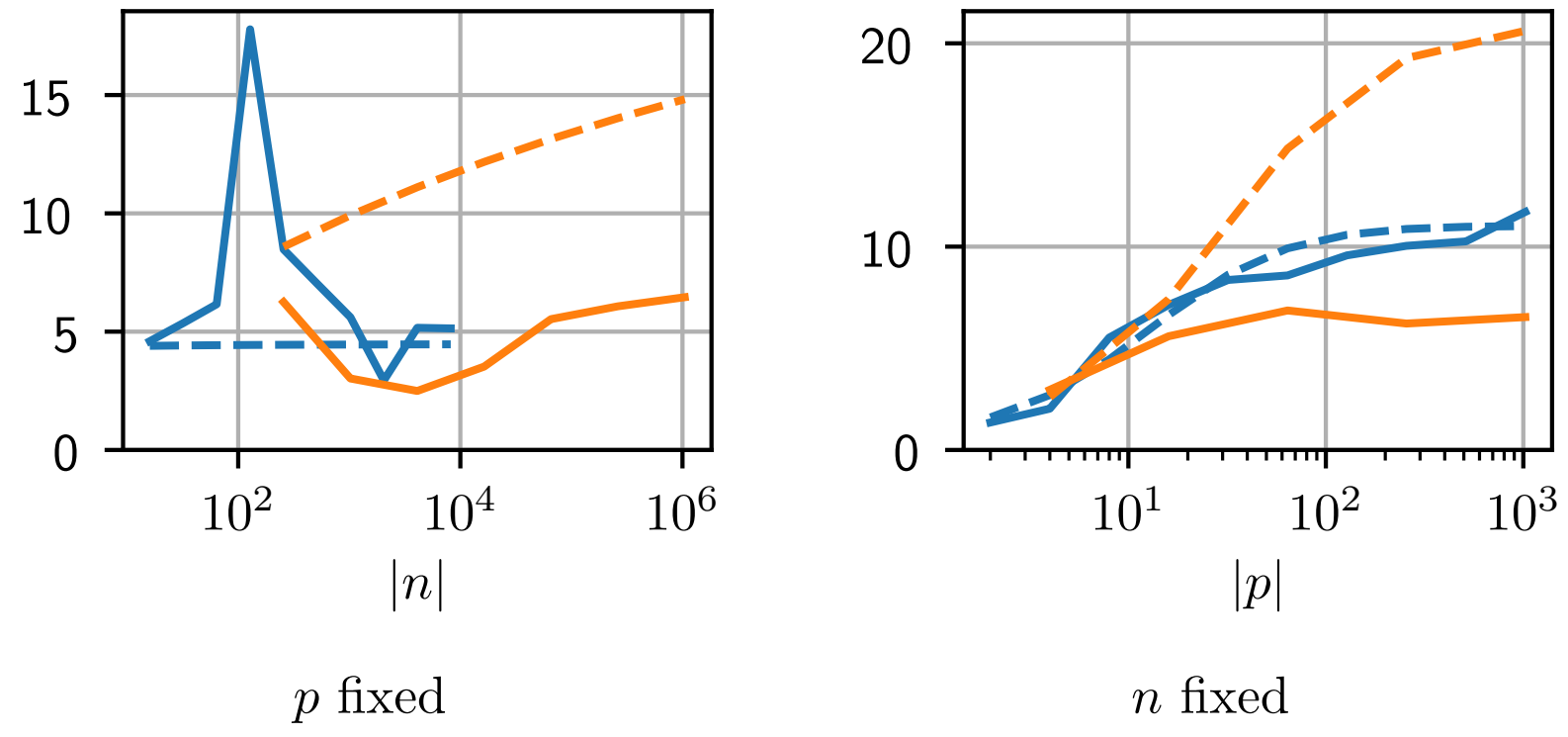

 \begin{tabular}{cc}
  \input{pfixed_speedup.pgf} &
  \input{nfixed_speedup.pgf} \\
  $p$ fixed & $n$ fixed
 \end{tabular}
 \caption{Speedups of the improved PCC algorithm compared to the
 circular convolution algorithm, for the periods and dimensions
 summarized in \cref{tab:summary}. The blue (resp. orange) curve is  for
 the 1D (resp. 2D) algorithm. The solid lines are the wallclock time
 speedups, whereas the dashed line is the one predicted by the cost
 estimates in \cref{tab:cost}.}
 \label{fig:speedup}
\end{figure}

We only include benchmarks for the DFT based algorithms as the
naïve algorithm becomes prohibitively expensive for relatively small
sizes. Wall clock timings were obtained on a MacBook Pro computer with
Apple M2 Max processor and 32 GB RAM which was running Mac OS 26.5 with
Julia 1.12.6. Moreover, both algorithms call extensively the DFT in the
FFTW 3.3.12 library. We performed two kinds of
experiments, one keeping $p$ fixed and varying $n$ and the other keeping
$n$ fixed and varying $p$.  The particular dimensions $n$ and periods
$p$ that we used are summarized in \cref{tab:summary} and all
calculations are carried with double precision complex floating point
numbers (16 bytes per entry).  The dimensions in 2D may appear small,
however a quick calculation reveals that simply storing the double
precision  complex  DFT of the kernel for the case where
$n=(2^{10},2^{10})$ and $p=(2^5,2^5)$ requires $2^{34} \approx 16$~GB of
memory. Average timings over 1000 trials are shown in \cref{fig:time}
while the speedups are in \cref{fig:speedup} and are compared to the
speedups predicted by dividing the application costs of the improved vs
convolution algorithms in \cref{tab:cost}.

Consistent with the cost estimates of \cref{sec:cost}, we observe the
improved PCC algorithm is faster than the circular convolution algorithm
in both 1D and 2D.  On closer inspection the speedups range from a
factor of $1$ to $17.5$. The predicted $p$ fixed 1D speedups seem to
underestimate the actual speedups, while they overestimate the actual 2D
speedups. For $n$ fixed, the predicted speedups for 1D closely match the
actual speedups, while in 2D the predicted ones are an overestimation,
see \cref{fig:speedup}. The comparisons between predicted vs actual
speedups should be taken with a grain of salt, since our cost estimates
were solely based on flop counts; and ignore memory access, allocation
costs and optimizations (processor or DFT library specific).

The improved DFT algorithm requires $ 3\times d$ nested loops, where
$d=1,2$ is the dimension ($3$ loops in 1D
and $6$ in 2D). In our implementation we chose to order the nested loops
so that the traversal of the kernel's DFT is done linearly.  Finding an
order to carry out the operations that minimizes the wall clock time may
depend on the dimensions of the signal, the period, the number of
threads, the number of cores and even the platform's memory hierarchy
(different caches etc.). Such an exploration is out of the scope of the
present work, but may benefit from an automatic tuning approach, similar
to that in e.g. FFTW \cite{Frigo:1998:FFT}. Both DFT PCC algorithms are
expected to benefit from a GPU implementation but this is also out of
the scope of this work.

\begin{table}
 \caption{Summary of the dimensions $n$ and periods $p$ we used in the
 numerical experiments.}
 \label{tab:summary}
 \begin{center}
 \begin{tabular}{|l|l|l|}
 \hline
  & $p$ & $n$\\\hline
  $p$ fixed 1D & $2^3$ &  $2^4,2^5,\ldots,2^{13}$\\
  $p$ fixed 2D & $(2^3,2^3)$ & $(2^4,2^4),\ldots,(2^{10},2^{10})$\\
  $n$ fixed 1D & $2^1,\ldots,2^{10}$ & $2^{10}$\\
  $n$ fixed 2D & $(2^1,2^1),\ldots,(2^5,2^5)$ &
  $(2^{10},2^{10})$\\\hline
 \end{tabular}
 \end{center}
\end{table}

\section{Linear algebra with PCCs}
\label{sec:blk}
The PCC operation being linear, it can be represented by a $|n| \times
|n|$ matrix that takes the vector representation $\vect(x) \in
\complex^{|n|}$ of the input signal to the vector representation of the
output signal $\vect(K *_p x) \in \complex^{|n|}$. Here $\vect(x)$ is
obtained by traversing linearly in memory the tensor $x$, assuming
column-major ordering. As we saw in \cref{rem:blk}, the matrix
representation $A \in \complex^{|n| \times |n|}$ of the PCC is of the
form:
\begin{equation}
 A = Q D_A Q^*,
  \label{eq:blkdiag}
\end{equation}
where $Q \in \complex^{|n| \times |n|} $ is a unitary matrix, obtained
by composition of DFTs and appropriate permutations, depending
\emph{only} on $n$ and $p$. The DFTs are \emph{normalized}, i.e. the
factor $1/\sqrt{|q|}$ is
included, so that $Q^*Q=I$ and \eqref{eq:blkdiag} holds as written. With
the unnormalized DFT \eqref{eq:dft}, the same composition $Q_0$ satisfies
$Q_0^*Q_0 = |q| I$, and \eqref{eq:blkdiag} reads $A = Q_0^* D_A
Q_0/|q|$. Moreover, $D_A \in \complex^{|n| \times
|n|}$ is a block diagonal matrix with blocks of size $|p| \times |p|$,
given by the unnormalized DFTs $\hat K$ of
\cref{eq:gconv1,eq:gconv2}.

Hence the matrix representations of PCCs (with same $n$ and $p$) form a
ring in $\complex^{|n| \times |n|}$ because block diagonal matrices with
same dimensions and block sizes form a ring. To be more precise, let
$A=QD_AQ^*$ and $B=QD_BQ^*$ be the block matrix representations of two
PCCs with same $n$ and $p$. Then their sum and product are  $A+B =
Q(D_A+D_B)Q^*$ and $AB = QD_AD_BQ^*$. Hence the block diagonal matrix
representation allows for calculation of the block diagonal matrix
representation of the product (composition) of two PCCs in $|q| |p|^3 =
|p|^2 |n|$ flops. Interestingly, the composition of PCCs is only
commutative in the case $|p|=1$, where the PCC reduces to a circular
convolution and $D_A$, $D_B$ are diagonal.

The adjoint of a PCC block diagonal matrix representation also represents
a PCC since $A^* = Q D_A^* Q^*$ and $D_A^*$ is block diagonal with
identical dimensions as $D_A$. It can be computed by conjugating and
appropriately permuting the dimensions in its tensor representation.
This is illustrated with functions \verb|pcc?d_adjoint| in
\texttt{PCC.jl}.

A PCC with a block diagonal matrix representation $A = QD_A Q^*$ is
invertible if and only if all the $|p| \times |p|$ diagonal blocks of
$D_A$ are invertible. Since it is a block diagonal matrix, the inverse
of $D_A$, if it exists, can be computed at a cost of $|q| |p|^3 = |p|^2
|n|$ flops. As a simple proof of concept, in \verb|pcc?d_solve| we show
how to solve linear systems with PCCs by solving a linear system for
each block, i.e. using the ``backslash'' operator in Julia.

The eigendecomposition of the matrix representation of a PCC is also
clear from its block diagonal representation. Indeed the
eigendecomposition of a Hermitian block diagonal matrix $D_A$ can be
found from the eigendecomposition of each individual diagonal block. The
eigenvectors of $D_A$ can be obtained by appropriate zero padding of the
eigenvectors of a particular diagonal block. Since $A$ is unitarily
equivalent to $D_A$, the eigenvectors of $A$ can be obtained from
those of $D_A$. In addition, the spectrum of $A$ is the union of the
spectra of the diagonal blocks of $D_A$. The Jordan canonical form of
$A$ can be obtained from that of $D_A$ in a similar fashion.
Generalized eigenvalues and eigenvectors for a linear pencil $(A,B)$
composed of matrices representing $p-$periodic convolutions can be
obtained from the matrix pencils of the block matrices in $D_A$ and $D_B$. 

Common matrix factorizations such as LU, Cholesky, QR, Schur and
Singular Value Decomposition can also be carried out block by block and
may be useful in applications, see e.g. \cite{Golub:2013:MC} for a
reference. For example, since the Frobenius norm and the induced matrix
$2-$norm (the largest singular value) are invariant through unitary
transformations, they can also be computed from the block diagonal
representation $D_A$. The Frobenius norm squared of $A$ is the sum of
the Frobenius norms of the diagonal blocks of $D_A$ and the induced
matrix $2-$norm of $A$ is the maximum of the union of singular values of
the diagonal blocks of $D_A$. Similarly, the $2-$norm condition number
of $A$ is given by the maximum to minimum ratio of the union of singular
values of the diagonal blocks of $D_A$ (this is implemented in
\verb|pcc?d_cond|). Orthogonal bases for the range
and nullspace of $A$, $A^*$ can be readily obtained from those of the
diagonal blocks $D_A$.  Finally, the rank of $A$ is simply the sum of
the ranks of each individual diagonal block in $D_A$.

\section{Deconvolution}
\label{sec:deconv}
We consider the problem of deconvolution, which consists of finding a
two dimensional signal $x$ from $y = K *_p x + w$, assuming the
$p-$periodic kernel $K$ is known. Here $w$ models noise in the
measurements. Let $A = Q D_A Q^*$ be the block diagonal representation
of the kernel $K$, obtained as in \cref{sec:blk}. Then the deconvolution
problem is that of finding $x$ from $y = A x + w$, where in a slight
abuse of notation, $x$, $y$ and $w$ are now vector representations. In general, the
matrix $A$ is not invertible. We could solve for $x$ by
minimizing the residual $\|y - A x\|_2^2$, but that could lead to
a spurious $x$ if the matrix $A$ is ill-conditioned. Instead, we use
Tikhonov regularized least squares, solving the optimization
problem
\begin{equation}
 \min_{x} \|A x - b\|_2^2 + \alpha^2 \|D x\|_2^2,
 \label{eq:tikho}
\end{equation}
where $\alpha>0$ is a regularization parameter and $D$ is an appropriate
matrix. The regularization parameter $\alpha$ serves to modulate the
importance we give to either minimizing the residual or the penalty term
$\|D x\|_2^2$, see e.g. \cite{Vogel:2002:CMI} for an introduction to
Tikhonov regularization. For the matrix $D$ we choose the one
corresponding to the discrete gradient $D = [ D_1; D_2]$, where $D_1$
(resp. $D_2$) are matrix representations of periodic differences in the first
(resp. second) dimension for a signal $x \in \complex^{n_1 \times
n_2}$. To be more precise we have
\begin{equation}
 D_1 x = x[1 \dv n_1,:] - x ~\text{and}~
 D_2 x = x[:,1 \dv n_2] - x,
 \label{eq:dop}
\end{equation}
recalling that array indices are periodic so the differences wrap
around in their respective dimensions. The choice of periodic difference
operators simplifies the problem to solve, as is observed in
\cite{Wang:2008:NAM}. Minimizing
\eqref{eq:tikho} can be done by solving for $x$ in the linear
system:
\begin{equation}
 (A^* A + \alpha^2 D_1^* D_1 + \alpha^2 D_2^* D_2) x = A^* y.
 \label{eq:normal}
\end{equation}
The key to solve \eqref{eq:normal} efficiently is to notice that all the
matrices involved are $p-$periodic convolutions. This is obviously the
case for $A^*$ (see \cref{sec:blk}). For $D_1$ and $D_2$, notice that
they are regular $p=(1,1)$ convolutions and therefore can also be
interpreted as $p-$periodic convolutions for any $p$ such that
$\mod{n_j}{p_j} = 0$, $j=1,2$. Indeed if $K$ is a two dimensional kernel
such that \eqref{eq:periodic:2d} holds for $p_1=p_2=1$, then
\eqref{eq:periodic:2d} also holds for any other valid $p$. We point out
that the matrix $D_1^* D_1 + D_2^* D_2$ is a graph Laplacian for the
connected graph formed by taking the Cartesian product of two cycle
graphs ($C_{n_1} \times C_{n_2}$), therefore it has a one-dimensional
nullspace spanned by the constant vector $1_{|n|}$, see e.g.
\cite{Chung:1997:SGT}. We conclude that the system matrix in
\eqref{eq:normal} is invertible provided the nullspace of $A$ does not
contain $1_{|n|}$. This is the case in our situation.  Indeed we see
from the singular values in \cref{fig:svd}(a) that $A$ has a trivial
nullspace containing only $0_{|n|}$.

\begin{figure}
 \begin{tabular}{ccc}
  \includegraphics[width=0.3\textwidth]{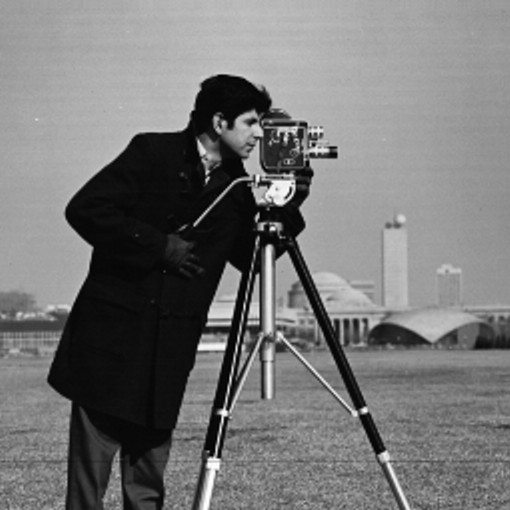} &
  \includegraphics[width=0.3\textwidth]{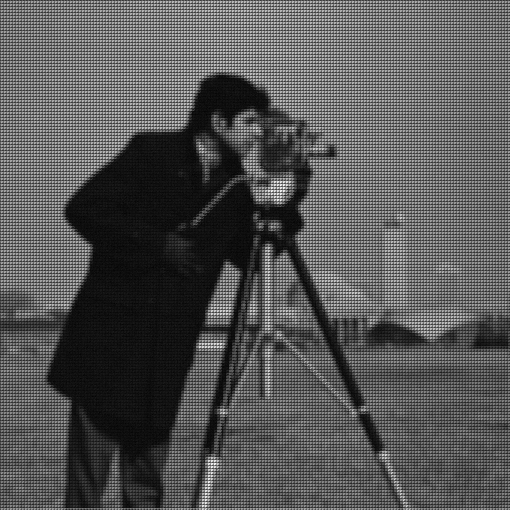} &
  \includegraphics[width=0.3\textwidth]{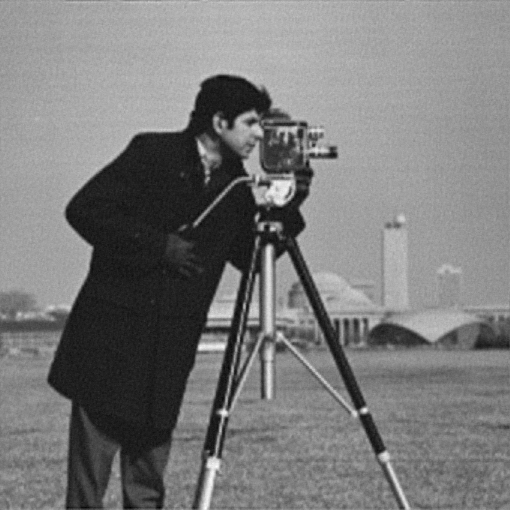} \\
  \includegraphics[width=0.3\textwidth]{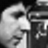} &
  \includegraphics[width=0.3\textwidth]{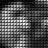} &
  \includegraphics[width=0.3\textwidth]{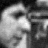} \\
  (a) & (b) & (c)
 \end{tabular}

 \caption{(a) True image $x_{\text{true}}$. (b) Right hand side
 $y_{\text{noisy}} = A x_{\text{true}} + $ noise. (c) Reconstructed~$x$. The
 top row shows the full $510\times 510$ pixel image while the bottom row is a $48 \times 48$ pixel inset.}
 \label{fig:deconv}
\end{figure}

\begin{figure}
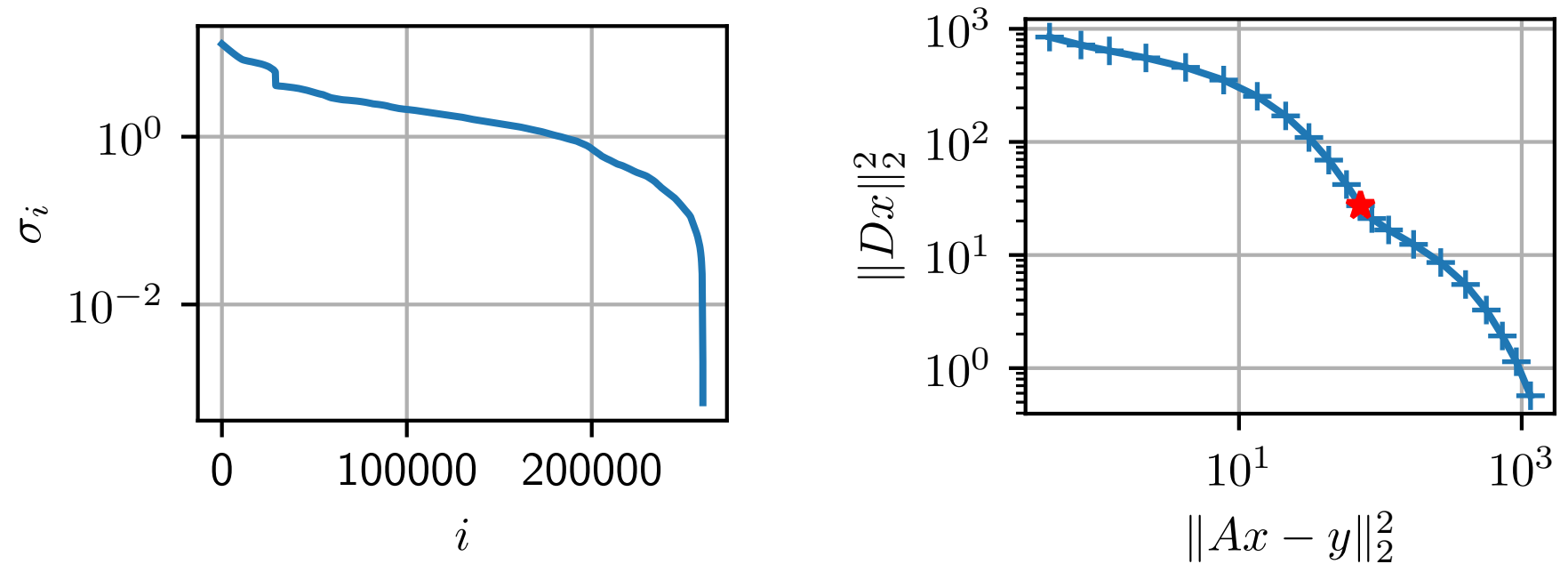

 \centering
 \begin{tabular}{cc}
  \input{deconv_svdvals.pgf} &
  \input{deconv_lcurve.pgf} \\
  (a) & (b)
 \end{tabular}
 \caption{(a) Singular values of system matrix $A$. Here $i$ is the
 singular value index. (b) L-curve
 generated for $\alpha = 2^{-10}, 2^{-9}, \ldots, 2^{10}$. The value 
 $\alpha=2$ we used in \cref{fig:deconv}(c) is in red.}
 \label{fig:svd}
\end{figure}

We illustrate this procedure with the same $(3,3)-$periodic kernel
presented in \cref{fig:pcc2d} but extended by zeros to $n=(510,510)$.
The original image $x_{\text{true}}$ is shown in \cref{fig:deconv}(a)
and is obtained from the well known \texttt{cameraman} test image by
removing the last two rows and columns so that $n$ is consistent with
$p=(3,3)$ and normalized to have values in $[0,1]$. The result $y = A x
+ w$ is displayed in \cref{fig:deconv}(b) and has a strong $3 \times 3$
pattern, as shown in the inset. The entries of the noise $w$ are
independent identically distributed normals with zero mean and standard
deviation equal to $\epsilon \|Ax\|_\infty$, where $\epsilon = 1/100$,
i.e. $1\%$ error relative to the largest entry of $A x$, in magnitude.
The singular values of $A$ can be calculated as explained in
\cref{sec:blk} and are given in \cref{fig:svd}(a). The $2-$norm
condition number for $A$ is $\kappa_2(A) \approx 18.8 \times 10^{3}$, so
the system is not singular but slightly ill-conditioned. The
reconstruction in \cref{fig:deconv}(c) is obtained with $\alpha = 2$
and has a $5.1\%$ relative $2-$norm error when compared to $x_{\text{true}}$.
The regularization parameter is determined via the L-curve method,
which is included in \cref{fig:svd}, see e.g. \cite{Vogel:2002:CMI} for
more details on this regularization parameter choice method. The idea is
to find an $\alpha$ that is a good balance between the residual $\|Ax -
b\|_2^2$ and the penalty $\|Dx\|_2^2$. This method is expensive because
it requires solving \eqref{eq:normal} for several values of $\alpha$.
However in the same platform used for \cref{sec:num}, solving
\eqref{eq:normal} takes about 100ms, without any
particular optimizations other than using the block matrix
representations introduced in \cref{sec:blk}.

\section{Dispersion relation for periodic chains}
\label{sec:greenper}
\begin{figure}
 \begin{center}
  \begin{tikzpicture}
  \foreach \i in {0,...,5} {
    \draw[thick, fill=black] (\i,0) circle (0.05);
  }
  \draw (0,0) -- (5,0);
  \draw[dashed] (5,0) -- (5.5,0);
  \draw[dashed] (-0.5,0) -- (0,0);
  \foreach \i in {0,...,2} {
   \foreach \j in {0,1} {
    \node[below] at (2*\i+\j,0) {$m_{\j}$};
   }
  }
  \end{tikzpicture}\\
  \input{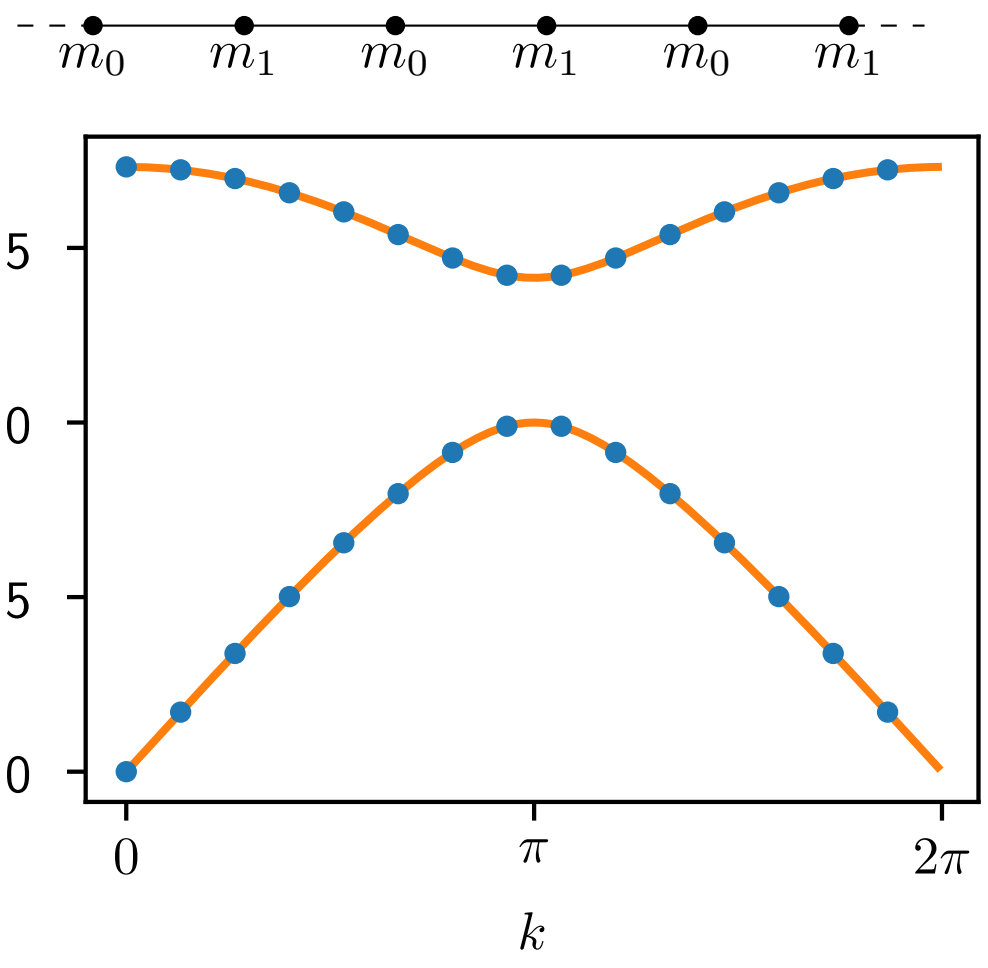}
 \end{center}
 \caption{Top: A periodic chain of $n=6$ masses with period $p=2$. The
 solid lines represent springs with spring constant $\kappa$ and the
 dashed lines is the spring between the first and last mass. Bottom: The
 continuous line is the true dispersion relation showing the acoustic
 $\omega_-(k)$ and optical $\omega_+(k)$ bands for $p=2$, $\kappa=1$,
 $m_0=1$ and $m_1=2$. The dots are the dispersion relation for $n=30$,
 calculated from the eigenvalues of the pencil $(K,M)$, indexed by the
 quasimomentum.}
 \label{fig:per}
\end{figure}

An elementary model from solid state physics is a chain of masses (or
atoms) that are arranged in a periodic fashion and linked with identical
springs, including between the last and first masses in the chain. We
use this model to point a connection between the block matrix
decomposition in \cref{sec:blk} and the Floquet-Bloch transform used to
analyze phenomena in periodic media, see e.g.
\cite{Berkolaiko:2013:IQG}. To fix ideas, consider a chain with $n$
masses that repeat with period $p$ (with $\mod{n}{p}=0$) and that are
linked by springs with spring constant $\kappa$. The case $n=6$ and
$p=2$ is given for instance in \cref{fig:per}(top). The vibrational
modes $u \neq 0$ of the chain satisfy 
\begin{equation}
 \kappa L u = \omega^2 M u,
 \label{eq:pencil}
\end{equation}
where $\omega$ is the angular frequency and $L$ is the
discrete graph Laplacian on the cycle graph $C_n$, i.e.
\begin{equation}
 L = 
 \begin{bmatrix}
  2 & -1 &&&-1\\
  -1 & 2 & -1\\
  & \ddots & \ddots & \ddots \\
  && -1 & 2 & -1\\
  -1&&& -1 & 2\\
 \end{bmatrix}
 ~\text{and}~
 M =
 \diag(m,m,\ldots,m),
\end{equation}
where $m = [m_0,m_1,\ldots,m_{p-1}]$.  Both $L$ and $M$ are
$p-$periodic convolutions since they satisfy \eqref{eq:periodic:1d}. In other
words, $L$ and $M$ commute with translations by an integer multiple of
$p$. Such translations form a group that is homomorphic to the cyclic
group $C_q$. The discrete Floquet-Bloch transform $\hat u \in
\complex^{n \times q}$ of $u$ can be defined as 
\begin{equation}
 \hat u [ i',i''] = \sum_{j''=0}^{q-1} u[i' + pi''] \exp[-\imath 2\pi i''
 j'' / q ],
\end{equation}
for $i' \in \irange{0,n-1}$, $i'' \in \irange{0,q-1}$ and $\imath =
\sqrt{-1}$.  Note that $\hat u$ is pseudoperiodic in the sense that for
any $i' \in \irange{0,n-1}$, $i'' \in \irange{0,q-1}$ and integer
$\ell$:
\begin{equation}
 \hat u[i' + \ell p , i'' ] = \exp[ \imath 2\pi i'' \ell / q] \hat u[i',i''].
\end{equation}
Hence $\hat u$ can be completely recovered from $\hat u[0 \dv p-1, :]$.
Upon closer examination we see that $\hat u [0 \dv p-1, :] = \cF_{*,q}
\reshape(u,p,q)$.  \Cref{thm:gconv1} shows that a $p-$periodic operator
$K$ is unitarily equivalent to a block diagonal matrix, with $p\times p$
blocks given by $\hat K [:,i'',:]$ for $i'' \in \irange{0,q-1}$.  The
\emph{dispersion relation} of $K$ is obtained by calculating the
eigenvalues of each block $\hat K[:,i'',:]$ and plotting them as a function
of the quasimomentum $2\pi i'' /q$, which plays the role of the
wavenumber. The vibration modes for the chain are given by the
eigenvalues of the matrix pencil $(\kappa L,M)$ as in \eqref{eq:pencil}.
Hence the dispersion relation for a period $p$ chain of atoms can be
obtained by plotting the eigenvalues of the real symmetric pencils
$(\kappa \hat L[:,i'',:],\hat M[:,i'',:])$ as a function of the quasimomentum
$2\pi i'' /q$. 

In the case $p=2$, the dispersion relation can be calculated by hand, as
can be seen in e.g. \cite{Ibach:2009:SSP}. We include here a sketch.
Let $\tilde u = \reshape(u,p,q)$ where $u \neq 0$ is a vibrational mode.
Then \eqref{eq:pencil} can be rewritten for $j\in \irange{0,q-1}$ as
\begin{equation}
\begin{aligned}
 \kappa ( -\tilde u[1,j] + 2\tilde u[0,j] - \tilde u[1,j-1]) &= m_0
 \omega^2 \tilde u[0,j], \\
 \kappa ( -\tilde u[0,j+1] + 2\tilde u[1,j] - \tilde u[0,j]) &= m_1
 \omega^2 \tilde u[1,j].
\end{aligned}
\end{equation}
Using the ansatz $\tilde u[:,j] = \exp[-\imath k j] A$ for
quasimomentum $k$ and some $A \in \complex^2$ we obtain
\begin{equation}
 \kappa
 \begin{bmatrix}
   2 & -(1+e^{\imath k})\\
   -(1+e^{-\imath k}) & 2
 \end{bmatrix}
 A
 = 
 \omega^2
 \begin{bmatrix}
  m_0 & \\ & m_1
 \end{bmatrix}
 A.
\end{equation}
Hence the natural vibration frequencies are given by
\begin{equation} 
\label{eq:dispersion} 
\omega^2_\pm(k) = \kappa (m_0^{-1} +m_1^{-1}) \pm \kappa \sqrt{(m_0^{-1}
+ m_1^{-1})^2  - \frac{4}{m_0m_1}\sin^2(\frac{k}{2})}. 
\end{equation} 
A comparison of the explicit expressions \eqref{eq:dispersion} with
those calculated block by block from \eqref{eq:pencil} is given in
\cref{fig:per}(bottom). Both results agree to machine precision. 

\section{Summary and conclusions}
\label{sec:conc}
We have presented two fast algorithms to compute periodic circular
convolutions (PCC) and which rely on the DFT. One of them requires
$|p|+1$ DFTs per application while the improved one requires only $2$
DFTs.
Depending on the size of the data and the period we observed
speedups of up to a factor of $17.5$ in the improved DFT PCC algorithm
versus the regular one. As discussed in \cref{sec:num}, the speedup
expected by counting flops is bounded by roughly $1+\log|n|$, but this
estimate falls short due in part to the 
importance of memory access.
The improved DFT PCC algorithm involves, aside from the DFT, a
tensor contraction operation that could be further optimized in terms of
memory access or by potentially exploiting the underlying matrix-vector
product structure.  Moreover, specialized algorithms for the DFT can
bring a factor of 2 improvement in both computational and memory usage
when the data is real \cite{Johnson:2026:IFP}.  One possible extension
of our approach would be to adapt the algorithms to take advantage of
these savings when both the kernel and the signal are real. 

The improved algorithm relies on a block diagonal representation of the
PCC.  This allows for a variety of linear algebra operations and
factorizations to be carried out in an efficient manner. In the
applications we show the inverse and generalized eigenvalue problems
were considered, but as mentioned in \cref{sec:blk}, other linear
algebra operations are possible.

\section*{Author Contributions} \textbf{Fernando González del Cueto}:
Conceptualization, Investigation, Software, Validation, Writing - Review
\& Editing. \textbf{Fernando Guevara Vasquez}: Conceptualization,
Investigation, Formal Analysis, Software, Visualization, Writing -
Original Draft.

\section*{Acknowledgements} The authors received funding from Lumos
Imaging Inc. FGV acknowledges funding from the National Science
Foundation DMS-2136198.

\bibliographystyle{siamplain}
\bibliography{pkc}

\end{document}